\documentclass[12pt,a4paper]{article}
\usepackage{amsmath,amssymb}
\usepackage{latexsym}
\usepackage{graphicx}
\usepackage{color}
\usepackage{multirow}
\usepackage{indentfirst}
\usepackage{subfigure}
\usepackage{mathrsfs}
\usepackage[pdftex]{hyperref}
 \usepackage[numbers,sort&compress]{natbib}

\newtheorem{exam}{\hspace{6mm}Example}[section]

\begin{document}
\baselineskip=2pc

\begin{center}
{\large \bf  A compact simple HWENO scheme with ADER time discretization for hyperbolic conservation laws I: structured meshes}
\footnote{The reasearch is supported partly by National Key R\&D Program of China (Grant Number 2022YFA1004500), National Natural Science Foundation of China (Grant Nos. 12101063, 12071392, 118872210, MCMS-I-0120G01).}
\end{center}

\centerline{
Dongmi Luo%
\footnote{Institute of Applied Physics and Computational Mathematics, Beijing 100088, China. E-mail: dongmiluo@stu.xmu.edu.cn.},
Shiyi Li%
\footnote{Institute of Applied Physics and Computational Mathematics, Beijing 100088, China. E-mail: lishiyi14@tsinghua.org.cn.},
Jianxian Qiu%
\footnote{School of Mathematical Sciences and Fujian Provincial
Key Laboratory of Mathematical Modeling and High-Performance
Scientific Computing, Xiamen University, Xiamen, Fujian 361005, China. E-mail: jxqiu@xmu.edu.cn.},
Jun Zhu%
\footnote{State Key Laboratory of Mechanics and Control of Mechanical Structures. Key Laboratory of Mathematical Modelling and High Performance Computing of Air Vehicles (NUAA), MIIT. Nanjing University of Aeronautics and Astronautics, Nanjing, Jiangsu 210016, P.R. China. E-mail: zhujun@nuaa.edu.cn.},
Yibing Chen*%
\footnote{Institute of Applied Physics and Computational Mathematics, Beijing 100088, China, 
E-mail: chen\_yibing@iapcm.ac.cn.}
}

\vspace{20pt}

\begin{abstract}
In this paper, a compact and high order ADER (Arbitrary high order using DERivatives) scheme using the simple HWENO method (ADER-SHWENO) is proposed for hyperbolic conservation laws. The newly-developed method 
employs the Lax-Wendroff procedure to convert time derivatives to spatial derivatives, which provides the time evolution of the variables at the cell interfaces. This information is required for the simple HWENO reconstructions, which take advantages of the simple WENO and the classic HWENO. Compared with the original Runge-Kutta HWENO method (RK-HWENO), the new method has two advantages. Firstly, RK-HWENO method must solve the additional equations for reconstructions, which is avoided for the new method. Secondly, the SHWENO reconstruction is performed once with one stencil and is different from the classic HWENO methods, in which both the function and its derivative values are reconstructed with two different stencils, respectively. Thus the new method is more efficient than the RK-HWENO method. Moreover, the new method is more compact than the existing ADER-WENO method. 
Besides, the new method makes the best use of the information in the ADER method.
Thus, the time evolution of the cell averages of the derivatives is simpler than that developed in the work [Li et. al., 447 (2021), 110661.].
 Numerical tests indicate that the new method can achieve high order for smooth solutions both in space and time, keep non-oscillatory at discontinuities.

\end{abstract}


\pagenumbering{arabic}

\newpage

\section{Introduction}
\label{sec1}
\setcounter{equation}{0}
\setcounter{figure}{0}
\setcounter{table}{0}
In this paper, we consider the numerical solutions of hyperbolic conservation laws. 
 The major difficulty in solving nonlinear  hyperbolic conservation laws is that the solution can develop discontinuities even if the initial condition is smooth. In the past decades, a large number of the high order finite volume methods have been developed to solve the equations, which can be divided for two groups. One is based on the cell averages, such as essentially non-oscillatory (ENO) \cite{harten1987}, weighted ENO (WENO) \cite{jiang1996,castro2011} and so on. However, the stencil used for the reconstructions is becoming wider with an increasing the order of accuracy. The other group is discontinuous Galerkin method (DG) \cite{cockburn1989jcp,cockburn1998,luo2019}, Hermite WENO (HWENO) \cite{qiu2003,qiu2005,zhu2009,luo2016,zhao2020} and so on, which are compact methods. Generally, these methods should solve the additional equations to obtain the internal freedom for the reconstructions, which is used to time evolution.
 
 Although these methods are arbitrarily high order in space. Most of them are discretized with Runge-Kutta (RK) methods in time. 
 For the multi-stage methods, the reconstruction needs to be performed several times per time step.
 To avoid the disadvantages of the multi-stage temporal discretization methods, lots of high order one-stage methods such GRP (the generalized Riemann problem) \cite{toro2009}, HGKS (high order gas kinetic scheme) \cite{li2021,liu2014,ren2015}, ADER (Arbitrary high order using DERivatives) approach of Toro et. al.  \cite{titarev2002,titarev2005,toro2006}, and other Lax-Wendroff procedure based methods \cite{qiu2003lw,qiu2005lw,qiu2007lw} are proposed. 
 
 In the past years, the gas-kinetic scheme (GKS) has been developed systematically \cite{li2021,liu2014,ren2015}. The high order GKS has been proposed without employing RK time discretization \cite{liu2014,ren2015}. 
 Moreover, a compact and efficient HGKS (CEHGKS) \cite{li2021} is proposed for the hyperbolic conservation laws recently, which is based on the framework of a one-stage efficient HGKS and can achieve high order accuracy both in space and time. 
The numerical examples \cite{li2021} show the effectiveness of CEHGKS. However, the technique is only suitable in the framework of HGKS.
 
Another popular one-stage temporal discretization method is ADER method, which was first put forward by Toro and collaborators for linear problems on Cartesian meshes \cite{toro2001}. Soon after that the ADER methods are extended to many nonlinear problems. ADER methods are combined with WENO methodology (ADER-WENO) to obtain the non-oscillatory of very high order of accuracy in space and time for nonlinear hyperbolic conservation laws \cite{titarev2002, titarev2005, toro2006, gu2017,dumbser20071,dumbser20072}.  Similar to the traditional RK-WENO, the stencil used in ADER-WENO reconstructions is becoming wider with an increasing the order of accuracy, which makes the method more complex on high dimension and unstructured meshes. To avoid the disadvantages of ADER-WENO methods, the ADER approaches are extended to nonlinear systems in the framework of DG methods (ADER-DG)\cite{dumbser2006,fambri2017}. However, ADER-DG methods also need to solve the additional equations for obtaining the internal freedom for reconstructions and time advancing, which increases the computational costs. 

Motivated by the CEHGKS and ADER methods, a compact and high order ADER scheme using the simple HWENO (SHWENO) method for hyperbolic conservation laws is proposed in this paper, which is denoted by ADER-SHWENO. The method is based on a fundamental assumption, i.e., all the variables of the equations are smooth within a space-time computational cell, which is also used for the work \cite{li2021} and a reason why the method \cite{li2021} can succeed indeed.
Thus the Lax-Wendroff procedure can be employed to convert time derivatives to spatial derivatives, which provides the time evolution of the variables at the cell interfaces. This information is required for the SHWENO reconstructions, which combine the ideas of the simple WENO \cite{zhu2016,zhu2017} and classic HWENO \cite{qiu2003,qiu2005,zhu2009}. Note that the classic HWENO \cite{qiu2003,qiu2005,zhu2009} must solve the additional equation for time evolution, which is avoided for the new method. Moreover, the SHWENO reconstruction is performed once with one stencil and is different from the classic HWENO methods \cite{qiu2003,qiu2005,zhu2009}, in which both the function and its derivative values are reconstructed with two different stencils, respectively. Thus ADER-SHWENO is more efficient than the classic HWENO methods. 
Moreover, the new method makes the best use of the information in the ADER method.  Therefore, the time evolution of the cell averages of the solutions is simpler than that in the work \cite{li2021}, where it includes the non-equilibrium and equilibrium parts, which need to be computed additionally and increase the computational cost. Besides, the approximation to leading term is $r$th order accuracy for the $r$th ADER-SHWENO method while $(2r-1)$th order accuracy for $r$th order ADER-WENO method. Therefore, the newly-developed ADER-SHWENO method is more compact than ADER-WENO method \cite{titarev2002}. Finally, an efficient and compact fifth order one-stage method is developed by combining ADER and SHWENO reconstructions. Numerical tests show that the new method can achieve high order for smooth solutions both in space and time, keep non-oscillatory at discontinuities. 

 An outline of the paper is given as follows. The ADER scheme based on the simple HWENO method is described in Section \ref{sec2} and Section \ref{sectwo} for one and two dimensional hyperbolic conservation laws in detail, respectively. In Section \ref{alg}, the algorithm of the method is presented briefly. One and two dimensional numerical examples are presented to demonstrate the accuracy and the capability of the ADER-SHWENO method in Section \ref{secnumexam}. In Section \ref{secconclusion}, conclusions are drawn. 

\section{The numerical scheme in one dimension}
\label{sec2}
\setcounter{equation}{0}
\setcounter{figure}{0}
\setcounter{table}{0}

In this section, we describe the ADER-SHWENO method for the numerical
solution of  hyperbolic conservation laws on a uniform mesh. 
We consider the conservation law in one dimension in the form 
\begin{equation}
\label{eq1}
\left\{
\begin{array}{l}
W_t + F(W)_x=0,\\
W(x,0)=W_0(x),
\end{array}
\right.
\end{equation}
where $W$ and $F(W)$ are either scalars or vectors.

Assume the computational domain $\Omega$ is divided into $N$ nonoverlapping cells $\{I_j=(x_{j-\frac{1}{2}},x_{j+\frac{1}{2}}),j=1,\cdots,N\}$, $\Delta x=x_{j+\frac{1}{2}}-x_{j-\frac{1}{2}}$, and $x_{j}=\frac{1}{2}(x_{j+\frac{1}{2}}+x_{j-\frac{1}{2}})$. Assume the computational time is $T$ and $0=t_0<t_1<\cdots <t_n<t_{n+1}<\cdots\leq T$. Integrating (\ref{eq1}) on $I_j\times(t_n,t_{n+1})$ yields
\[
\Bar W(x_j,t_{n+1})=\Bar W(x_j,t_{n})-\frac{\Delta t}{\Delta x}(\frac{1}{\Delta t}\int_{t_n}^{t_{n+1}}F(x_{j+\frac{1}{2}},t)dt-\frac{1}{\Delta t}\int_{t_n}^{t_{n+1}}F(x_{j-\frac{1}{2}},t)dt),
\]
where $\Bar W(x_j,t_{n})=\frac{1}{\Delta x}\int_{I_j}W(x,t_n)dx$ is the cell average in cell $I_j$ at time $t_n$, and $\Delta t=t_{n+1}-t_n$.
 Thus, the  finite volume method for the equation is defined as
\begin{align}
\label{cons}
\Bar W_j^{n+1}=\Bar W_j^{n} -\frac{\Delta t}{\Delta x}(\Hat F_{j+\frac{1}{2}}-\Hat F_{j-\frac{1}{2}}),
\end{align}
where $\Bar W_j^n$ is a high order approximation to $\Bar W(x_j,t_n)$, and $\Hat F_{j+\frac{1}{2}}\approx \frac{1}{\Delta t}\int_{t_n}^{t_{n+1}}F(x_{j+\frac{1}{2}},t)dt$ is the numerical flux, which is also a high order approximation to the physical flux. For the traditional finite volume methods, an explicit and the third order TVD Runge-Kutta scheme \cite{shu1988} is used for temporal discretization. In order to obtain the arbitrary order of accuracy in both time and space, the ADER approach \cite{titarev2002} is employed in this paper. The ADER approach mainly contains three steps \cite{titarev2002}: 

1. first reconstruct the high order point-wise values from the cell averages; 

2. after the reconstruction, solve the generalize Riemann problem at the cell interface;

 3. finally evaluate the numerical flux in the conservative scheme (\ref{cons}).
 
%
%

\subsection{The generalized Riemann problem}
\label{grp}
In this section, we describe the generalized Riemann problem. We assume the point-wise values of the conservative variables have been reconstructed, which are represented by the vector polynomials $p_j(x), j=1,\cdots,N$. Then at each cell interface, the generalized Riemann problem (GRP) with the reconstruction polynomials is proposed \cite{titarev2005,titarev2002}:
\begin{equation}
\label{grp1d}
\begin{array}{ll}
\partial_t W+\partial_xF(W) =0 \\
W(x,0)=
\left
\{
\begin{array}{ll}
W_L(x)=p_j(x),&\; \;\; x<x_{j+\frac{1}{2}},\\
W_R(x)=p_{j+1}(x),&\;\;\;  x>x_{j+\frac{1}{2}},
\end{array}
\right.
\end{array}
\end{equation}
which is different from the conventional piece-wise constant data Riemann problem. Now the solution no longer contains regions of constant values and characteristics are curved lines.
Thus, we need to find an approximate solution for the interface state $W(x_{j+\frac{1}{2}},\tau)$, where $\tau$ is local time $\tau=t-t_n$ using the method developed in \cite{toro2009}. The approximate solution $W(x_{j+\frac{1}{2}},\tau)$ can be evaluated by a Taylor expansion of the interface state in time
\begin{equation}
\label{tay}
W(x_{j+\frac{1}{2}},\tau)\approx W(x_{j+\frac{1}{2}},0^+)+\sum\limits_{k=1}^{4}[\partial_t^{(k)}W(x_{j+\frac{1}{2}},0^+)]\frac{\tau^k}{k!},
\end{equation}
where
\[
\partial_t^{(k)}W(x,t)=\frac{\partial^k}{\partial t^k}W(x,t),\quad 0^+\equiv \lim_{\tau\rightarrow 0^+}\tau.
\]
(\ref{tay}) implies the solution is smooth in time.
From (\ref{tay}), one can observe that the approximate solution contains two parts: a leading term $W(x_{j+\frac{1}{2}},0^+)$ and the high order terms with coefficients determined by the temporal derivatives $\partial_t^{(k)}W(x_{j+\frac{1}{2}},0^+)$. Thus, the values of a leading term and high order terms at the cell interfaces need to be evaluated, which are described in the following sections.
\subsubsection{The leading term}
The leading term $W(x_{j+\frac{1}{2}},0^+)$ accounts for the interaction of the boundary extrapolated values $W_L(x_{j+\frac{1}{2}})$ and $W_R(x_{j+\frac{1}{2}})$, which are taken to be the Godunov state of the conventional Riemann problem: 
\begin{equation}
\begin{array}{ll}
\partial_t W+\partial_xF(W) =0 \\
W(x,0)=
\left
\{
\begin{array}{ll}
W_L(x_{j+\frac{1}{2}}), &\; \;\; x<x_{j+\frac{1}{2}},\\
W_R(x_{j+\frac{1}{2}}),&\;\;\;  x>x_{j+\frac{1}{2}}.
\end{array}
\right.
\end{array}
\label{leadingterm}
\end{equation}
The leading term is obtained by solving the above equation. In this paper, $W_L(x_{j+\frac{1}{2}})$ and $W_R(x_{j+\frac{1}{2}})$ are directly obtained by the SHWENO reconstruction in Sec. \ref{sechweno}. For scalar cases, we solve the Riemann problem exactly. For the Euler systems, the HLLC Riemann solver \cite{toro2009} is adopted, which contains all waves in the Riemann problem solution.

\subsubsection{The high order term}
To compute the high order terms in (\ref{tay}) we need to compute the coefficients, which are the partial derivatives $\partial_t^{(k)}W(x_{j+\frac{1}{2}},0^+)$ at $x=x_{j+\frac{1}{2}}, t=0$. Thus, the high order terms are evaluated in two steps. 

First, we express all time derivatives as functions of space derivatives by the Lax-Wendroff  procedure. For systems the procedure yields the following expressions:
\begin{align}
\label{ck}
\partial_tW&=-\frac{\partial F}{\partial W}\partial_xW,\notag\\
\partial_{tx}W&=-(\frac{\partial^2 F}{\partial W^2}\partial_xW)\partial_xW-\frac{\partial F}{\partial W}\partial_{xx}W,\\
\partial_{tt}W&=-(\frac{\partial^2 F}{\partial W^2}\partial_tW)\partial_xW-\frac{\partial F}{\partial W}\partial_{xt}W,\notag
\end{align}
and so on. In practice, it is more convenient to implement (\ref{ck}) in componentwise manner rather than in the matrix form.

Then differentiating the governing equation (\ref{eq1}) with respect to $x$, we obtain
\begin{equation}
\partial_tW^{(k)}+A(W)\partial_xW^{(k)}=H(W,W^{(1)},W^{(2)},\cdots,W^{(k-1)}),
\end{equation}
where $W^{(k)}=\frac{\partial^k}{\partial x^k}W, 1\leq k \leq 4$, $A(W)$ is the Jacobian matrix, $H$ is a nonlinear source term depending on derivatives of lower order as well as $W(x, t)$ itself. However, it is not easy to solve these nonlinear inhomogeneous problems. For simplicity, we neglect the source terms $H$ which comes into effect for $\tau>0$ only. Additionally, we linearize the equation around the leading term $W(x_{j+\frac{1}{2}},0^+)$ of the time expansion (\ref{tay}) and replace the piece-wise polynomial initial data by left and right boundary extrapolated values of spatial derivatives at $x_{j+\frac{1}{2}}$. Then high order terms are obtained by solving the following linearized Riemann problem \cite{titarev2002, toro2006}:
\begin{equation}
\label{lineq}
\begin{array}{ll}
\partial_t W^{(k)}+A_{j+\frac{1}{2}}\partial_xW^{(k)} =0 \\
W^{(k)}(x,0)=
\left
\{
\begin{array}{ll}
\partial ^{(k)}_xW_L(x_{j+\frac{1}{2}}), &\; \;\; x<x_{j+\frac{1}{2}},\\
\partial ^{(k)}_xW_R(x_{j+\frac{1}{2}}),&\;\;\;  x>x_{j+\frac{1}{2}},
\end{array}
\right.
\end{array}
\end{equation}
where $A_{j+\frac{1}{2}}=F'(W(x_{j+\frac{1}{2}},0^+))$ is the Jacobian matrix. 

Note that the coefficient matrix $A_{j+\frac{1}{2}}$ is the same for all spatial derivatives $W^{(k)}$ and is evaluated only once, using the leading term of the expansion. The initial condition for (\ref{lineq}) is found directly by differentiating the given SHWENO 
reconstruction polynomial with respect to $x$. In the paper, we use the same linear weights and smoothness indicators for the function and for all derivatives. Then the linear Riemann problem \eqref{lineq} can be solved directly. 

\subsection{The evaluation of the numerical flux}
Finally, having computed all spatial derivatives we form the Tayor expansion (\ref{tay}). There exist two options to evaluate the numerical flux \cite{titarev2005}. In the paper, the state-expansion ADER \cite{titarev2002} is employed, in which the approximate state (\ref{tay}) is inserted into the definition of the numerical flux (\ref{cons}). To evaluate the numerical flux an appropriate Gaussian rule is used:
\begin{equation}
\label{flux}
\Hat F_{j+\frac{1}{2}}=\sum\limits_{\alpha=0}^{K_{\alpha}}F(W(x_{j+\frac{1}{2}},\lambda_{\alpha}\Delta t))w_{\alpha},
\end{equation}
where $\lambda_{\alpha}$ and $w_{\alpha}$ are properly scaled nodes and weights of the rule, and $K_{\alpha}$ is the number of the nodes.

\subsection{SHWENO reconstruction in one dimension}
\label{sechweno}
In the previous sections, the GRP and evaluation of the flux are presented. Recall that we assume the polynomials are reconstructed. In this subsection, the simple HWENO reconstruction is given. For the traditional ADER methods, both WENO and DG methods are applied for reconstructions. For the classic WENO method \cite{jiang1996}, the approximation to the leading term is $(2r-1)$th order accuracy for $r$th order ADER method. Thus, the stencil used in reconstructions is becoming wider with an increasing order of accuracy. On the other hand, for DG method \cite{cockburn1989jcp,cockburn1998}, the additional equations need to be solved to obtain the internal freedom for reconstruction and time advancing, which increases the computational costs. To avoid the disadvantages of WENO and DG method, the simple HWENO method \cite{li2021} is employed in this paper. The reconstruction method combine the ideas of the simple WENO and the classic HWENO method. Thus the linear weights can be chosen arbitrary positive numbers except that its sum equal to one and the cell averages as well as the cell averages of the derivatives are needed for reconstruction.  The procedure is the same as the work \cite{li2021} and is omitted here. In this work, we set the linear weights as follows: $\gamma_1=0.994, \gamma_2=0.003$ and $\gamma_3=0.003$.

\subsection{The evaluation of the cell average of $W_x$}
Note that in Sec. \ref{sechweno} the cell averages of $W_x$ is needed for the SHWENO reconstruction in the next time step, which is evaluated in the following. Assume 
\begin{equation}
\label{vv}
\Bar V_j=\frac{1}{\Delta x}\int_{I_j} W_x(x,t)dx,
\end{equation}
which is the cell average of the derivatives of the solution. From (\ref{vv}), at time $t_{n+1}$ we have
\begin{equation}
\label{derave}
\Bar V^{n+1}_j=\frac{1}{\Delta x}\int_{I_j}W_x(x,t_{n+1})dx=\frac{W(x_{j+\frac{1}{2}},t_{n+1})-W(x_{j-\frac{1}{2}},t_{n+1})}{\Delta x}.
\end{equation}
In the classic HWENO method \cite{qiu2003}, an additional equation with respect to the cell average of $W_x$ is solved for time evolution. To avoid the disadvantages of HWENO methods,
 Li et al. \cite{li2021} employed the equation (\ref{derave}) and evaluated the value of $W(x_{j+\frac{1}{2}},t_{n+1})$ in the framework of GKS, which includes the non-equilibrium and equilibrium states. Thus these terms should be evaluated in order to calculate $W(x_{j+\frac{1}{2}},t_{n+1})$, which increases the complexity of the method developed in \cite{li2021}. 


Recall the ADER method described in the previous section. From \eqref{flux}, one can see that the numerical solutions at the different times between $t_n$ and $t_{n+1}$ should be computed for the numerical flux. In the ADER method, the solutions at different times are evaluated by a Taylor expansion of the interface state in time, i.e. \eqref{tay}. Motivated by that, a fifth order approximation of $W(x_{j+\frac{1}{2}},t_{n+1})$ can be evaluated by
\begin{equation}
\label{tay1d}
W(x_{j+\frac{1}{2}},t_n+\Delta t)\approx W(x_{j+\frac{1}{2}},t_n^+)+\sum\limits_{k=1}^{4}[\partial_t^{(k)}W(x_{j+\frac{1}{2}},t_n^+)]\frac{(\Delta t)^k}{k!},
\end{equation}
which implies the solution in time is smooth.

In fact, the leading term $W(x_{j+\frac{1}{2}},t_n^+)$ and high order term $\partial_t^{(k)}W(x_{j+\frac{1}{2}},t_n^+)$ have been evaluated in the ADER method, which can be used directly in \eqref{tay1d}. Therefore, the evaluation of $W(x_{j+\frac{1}{2}},t_{n+1})$ in our new method almost does not introduce the additional computational cost and is simpler than that in the work \cite{li2021}. Finally, we can obtain the cell averages of $W_x$ at the next time using \eqref{derave}.

{\bf Remark 1:} From the procedure above, one can see that we only need to take $\tau=\Delta t$ in (\ref{tay}) to get \eqref{tay1d}  for updating the cell average of $W_x$. In addition, the leading term and high order term are evaluated in ADER method. Therefore, we almost directly use trivial work to update the cell average of $W_x$ using (\ref{derave}). This is the key point of our method.

{\bf Remark 2:} In the practical implementations, we can use the Gauss-Lobatto quadrature rule in \eqref{flux}, which includes the values at both ends, i.e. $t_n$ and $t_{n+1}$. Thus, we can further reduce the computational cost.

\section{The numerical scheme in two dimensions}
\label{sectwo}
\setcounter{equation}{0}
\setcounter{figure}{0}
\setcounter{table}{0}
In this section, we describe the ADER-SHWENO method for two dimensional problems. We consider the form
\begin{equation}
\label{eq2}
\left\{
\begin{array}{l}
W_t + F(W)_x + G(W)_y=0,\\
W(x,y,0)=W_0(x,y),
\end{array}
\right.
\end{equation}
where $W$, $F(W)$ and $G(W)$ are either scalars or vectors.
We use a rectangle mesh of cell size $\Delta x$ and $\Delta y$ in $x$ and $y$ directions, respectively. We denote the cells by
\[
I_{ij}=(x_{i-\frac{1}{2}},x_{i+\frac{1}{2}})\times (y_{j-\frac{1}{2}},y_{j+\frac{1}{2}}),
\]
where $
x_{i+\frac{1}{2}}=\frac{1}{2}(x_i + x_{i+1}), y_{j+\frac{1}{2}}=\frac{1}{2}(y_j + y_{j+1}).
$

Integrating \eqref{eq2} on $I_{ij}\times (t_n,t_{n+1})$, we obtain
\begin{align*}
\Bar W(x_i,y_j,t_{n+1})=\Bar W(x_i,y_j,t_{n})-&\frac{\Delta t}{\Delta x\Delta y}[\frac{1}{\Delta t}\int_{t_n}^{t_{n+1}}\int_{y_{j-\frac{1}{2}}}^{y_{j+\frac{1}{2}}}F(W(x_{i+\frac{1}{2}},y,t))-F(W(x_{i-\frac{1}{2}},y,t))dydt]-\\
&\frac{\Delta t}{\Delta x\Delta y}[\frac{1}{\Delta t}\int_{t_n}^{t_{n+1}}\int_{x_{i-\frac{1}{2}}}^{x_{i+\frac{1}{2}}}G(W(x,y_{j+\frac{1}{2}},t))-G(W(x,y_{j-\frac{1}{2}},t))dxdt],
\end{align*}
where $\Bar W(x_i,y_j,t_{n})=\frac{1}{\Delta x\Delta y}\int_{I_{ij}}W(x,y,t_n)dxdy$. Then the finite volume method in two dimensions is given by
\begin{equation}
\Bar W_{ij}^{n+1}=\Bar W_{ij}^{n}-\frac{\Delta t}{\Delta x\Delta y}(\Hat F_{i+\frac{1}{2},j}-\Hat F_{i-\frac{1}{2},j})-\frac{\Delta t}{\Delta x\Delta y}(\Hat G_{i,j+\frac{1}{2}}-\Hat G_{i,j-\frac{1}{2}})
\end{equation}
where $\Hat F_{i+\frac{1}{2},j}\approx \frac{1}{\Delta t}\int_{t_n}^{t_{n+1}}\int_{y_{j-\frac{1}{2}}}^{y_{j+\frac{1}{2}}}F(W(x_{i+\frac{1}{2}},y,t))dydt, \Hat G_{i,j+\frac{1}{2}}\approx \frac{1}{\Delta t}\int_{t_n}^{t_{n+1}}\int_{x_{i-\frac{1}{2}}}^{x_{i+\frac{1}{2}}}G(W(x,y_{j+\frac{1}{2}},t))dxdt$. One can observe that this leads to equations including the line integrals on the cell $I_{ij}$, which can be computed by the Gaussian quadrature rule
\begin{align}
\label{fx}
\Hat F_{i+\frac{1}{2},j}&=\sum\limits_{G_y}\Delta y w_{G_y}(\frac{1}{\Delta t}\int_{t_n}^{t_{n+1}}F(W(x_{i+\frac{1}{2}},y_{G_y},t))dt),\\
\label{gy}
 \Hat G_{i,j+\frac{1}{2}}&=\sum\limits_{G_x}\Delta x w_{G_x}(\frac{1}{\Delta t}\int_{t_n}^{t_{n+1}}G(W(x_{G_x},y_{i+\frac{1}{2}},t))dt),
\end{align}
where $x_{G_x}$ and $y_{G_y}$ are the Gaussian points on $(x_{i-\frac{1}{2}},x_{i+\frac{1}{2}})$ and $(y_{j-\frac{1}{2}},y_{j+\frac{1}{2}})$, respectively. $w_{G_x}$ and $w_{G_y}$ are the corresponding weights. As one dimensional case, we should evaluate the numerical fluxes \eqref{fx} and \eqref{gy}. Here we concentrate on $\Hat F_{i+\frac{1}{2},j}$ and the other one $ \Hat G_{i,j+\frac{1}{2}}$ is obtained in a similar manner.

\subsection{The generalized Riemann problem in two dimensions}
After the reconstruction is carried out for each Gaussian integration point $(x_{i+\frac{1}{2}},y_{G_y})$ at the cell interface, the generalized Riemann problem (\ref{grp1d}) can be proposed along the $x$-direction. Then we can obtain a high order approximation to $W(x_{i+\frac{1}{2}},y_{G_y},\tau)$. All steps of the solution procedure are the same as in the one dimensional case. The Taylor series expansion in time is written as following:
\begin{align}
\label{taylor2d}
W(x_{i+\frac{1}{2}},y_{G_y},\tau)=W(x_{i+\frac{1}{2}},y_{G_y},0^+)+\sum\limits_{k=1}^{4}[\partial_t^{(k)}W(x_{j+\frac{1}{2}},y_{G_y},0^+)]\frac{\tau^k}{k!}.
\end{align}
The leading term $W(x_{i+\frac{1}{2}},y_{G_y},0^+)$ is the Godunov state of the conventional Riemann problem:
\begin{equation}
\begin{array}{ll}
\partial_t W+\partial_xF(W) =0 \\
W(x,0)=
\left
\{
\begin{array}{ll}
W_L(x_{i+\frac{1}{2}},y_{G_y}), &\; \;\; x<x_{i+\frac{1}{2}},\\
W_R(x_{i+\frac{1}{2}},y_{G_y}),&\;\;\;  x>x_{i+\frac{1}{2}},
\end{array}
\right.
\end{array}
\label{leadingterm2d}
\end{equation}
which can be solved by an exact or HLLC Riemann solver. To evaluate the high order term we first use the Lax-Wendroff procedure to express the time derivatives as functions of space derivatives. Then we can obtain
\begin{align*}
\partial_t W &= -(\frac{\partial F}{\partial W})\partial_x W - (\frac{\partial G}{\partial W})\partial_y W,\\
\partial_{tx}W &= -(\frac{\partial^2 F}{\partial W^2}\partial_x W)(\partial_x W)  -(\frac{\partial F}{\partial W})\partial_{xx} W-(\frac{\partial^2 G}{\partial W^2}\partial_x W)(\partial_y W)-(\frac{\partial G}{\partial W})\partial_{xy} W,\\
\partial_{ty}W &= -(\frac{\partial^2F}{\partial W^2}\partial_y W)(\partial_x W) -(\frac{\partial F}{\partial W})\partial_{xy} W-(\frac{\partial^2 G}{\partial W^2}\partial_y W)(\partial_y W)-(\frac{\partial G}{\partial W})\partial_{yy} W,\\
\partial_{tt}W &= -(\frac{\partial^2 F}{\partial W^2}\partial_t W)(\partial_x W)  -(\frac{\partial F}{\partial W})\partial_{tx} W-(\frac{\partial^2 G}{\partial W^2}\partial_t W)(\partial_y W)-(\frac{\partial G}{\partial W})\partial_{ty} W,
\end{align*}
and so on. These equations can be also used for systems. However, $\frac{\partial F}{\partial W}$ is a matrix and $\frac{\partial^2 F}{\partial W^2}$ is a three-dimensional tensor, etc. After getting the expression, the high order terms are obtained by solving the linearized Riemann problem:
\begin{equation}
\label{lineq2d}
\begin{array}{ll}
\partial_t W^{(m+n)}+A_{i+\frac{1}{2},j}\partial_xW^{(m+n)} =0 \\
W^{(m+n)}(x,y_{G_y},0)=
\left
\{
\begin{array}{ll}
\partial ^{(m+n)}_{x^my^n}W_L(x_{i+\frac{1}{2}},y_{G_y}), &\; \;\; x<x_{i+\frac{1}{2}},\\
\partial ^{(m+n)}_{x^my^n}W_R(x_{i+\frac{1}{2}},y_{G_y}),&\;\;\;  x>x_{i+\frac{1}{2}},
\end{array}
\right.
\end{array}
\end{equation} 
where $W^{(m+n)}=\frac{\partial^{m+n}}{\partial x^m\partial y^n}W=\partial ^{(m+n)}_{x^my^n}W$, $1\leq m+n \leq 4$ and $A_{i+\frac{1}{2},j}=F'(W(x_{i+\frac{1}{2}},y_{G_y},0^+))$. Then the Taylor expansion \eqref{taylor2d} can be formed for the interface state at the Gaussian point $(x_{i+\frac{1}{2}},y_{G_y})$. Finally the numerical flux can be evaluated by
\begin{align}
\label{fx2d}
\Hat F_{i+\frac{1}{2},j}&=\sum\limits_{G_y}\Delta y w_{G_y}\sum\limits_{\alpha=0}^{K_{\alpha}}F(W(x_{i+\frac{1}{2}},y_{G_y},\lambda_{\alpha}\Delta t))w_{\alpha},
\end{align}

\subsection{SHWENO reconstruction in two dimensions}
\label{sechweno2d}
For SHWENO reconstruction on Cartesian meshes, one can employ either a direct two dimensional procedure or a dimension-by-dimension strategy \cite{luo2016}. In this paper, the dimension-by-dimension strategy is adopted. To perform the reconstruction,  the cell averages of $W, W_x, W_y$ and $W_{xy}$ are denoted by
\begin{equation}
\begin{array}{llll}
\Delta x\Delta y\Bar W_{ij}&=\int_{I_{ij}}Wdxdy,\\ 
\Delta x\Delta y\Bar V_{ij}&=\int_{I_{ij}}W_xdxdy,\\
\Delta x\Delta y\Bar Y_{ij}&=\int_{I_{ij}}W_ydxdy,\\
\Delta x\Delta y\Bar Z_{ij}&=\int_{I_{ij}}W_{xy}dxdy.
\end{array}
\end{equation}
The reconstruction is as follows. First, we perform two $y$-direction reconstructions as Sec. \ref{sechweno}, i.e.,
\begin{align*}
\{\Bar W_{mn}, \Bar Y_{mn} \}&\rightarrow \Bar W_{i+l,j}(y_{G_y})\approx \frac{1}{\Delta x}\int_{I_{i+l},j}W(x,y_{G_y})dx,\\
\{\Bar V_{mn}, \Bar Z_{mn} \}&\rightarrow \Bar W_{x,i+l,j}(y_{G_y})\approx \frac{1}{\Delta x}\int_{I_{i+l},j}W_x(x,y_{G_y})dx.
\end{align*}
Then we use $\Bar W_{i+l,j}(y_{G_y})$ and $\Bar W_{x,i+l,j}(y_{G_y})$ to perform $x$-direction reconstruction to get an approximation to $W(x_{i+\frac{1}{2}},y_{G_y})$, i.e.,
\begin{align*}
\{\Bar W_{mn}(y_{G_y}), \Bar W_{x,mn}(y_{G_y}) \} \rightarrow \tilde W(x_{i+\frac{1}{2}},y_{G_y})\approx W(x_{i+\frac{1}{2}},y_{G_y}).
\end{align*}
The reconstruction of $W(x_{G_x},y_{j+\frac{1}{2}})$ is performed in a similar way.

\subsection{The evaluation of the cell averages of $W_x$, $W_y$ and $W_{xy}$}
As in one dimensional case, the cell averages of $W_x$, $W_y$ and $W_{xy}$ need to be evaluated at the next time, which are given by
\begin{equation}
\label{derave2d}
\begin{array}{llll}
\Bar V_{ij}^{n+1}=&\frac{1}{\Delta x}\sum\limits_{G_y}w_{G_y}(W(x_{i+\frac{1}{2}},y_{G_y},t_{n+1})-W(x_{i-\frac{1}{2}},y_{G_y},t_{n+1})),\\
\Bar Y_{ij}^{n+1}=&\frac{1}{\Delta y}\sum\limits_{G_x}w_{G_x}(W(x_{G_x},y_{j+\frac{1}{2}},t_{n+1})-W(x_{G_x},y_{j-\frac{1}{2}},t_{n+1})),\\
\Bar Z_{ij}^{n+1}=&\frac{1}{\Delta x\Delta y}(W(x_{i+\frac{1}{2}},y_{j+\frac{1}{2}},t_{n+1})-W(x_{i-\frac{1}{2}},y_{j+\frac{1}{2}},t_{n+1})-\\
&\qquad\; (W(x_{i+\frac{1}{2}},y_{j-\frac{1}{2}},t_{n+1})-W(x_{i-\frac{1}{2}},y_{j-\frac{1}{2}},t_{n+1}))).
\end{array}
\end{equation}
 Similar to one dimensional case, Li et al. \cite{li2021} computed the values at next time using the non-equilibrium and equilibrium states.
 However, what is different from one dimensional case is that the computational costs are significantly increased in two dimensions since more terms at the next time should be evaluated for the numerical flux. Thus, the method developed in \cite{li2021} becomes more complicated in two dimensions.

%

In this paper, to get $W(x_{i+\frac{1}{2}},y_{G_y},t_{n+1})$, $\tau$ is taken as $\Delta t$ in (\ref{taylor2d}) similar to the one dimensional case, which leads
\begin{align}
\label{taylor2dtime}
W(x_{i+\frac{1}{2}},y_{G_y},t_n+\Delta t)=W(x_{i+\frac{1}{2}},y_{G_y},t_n^+)+\sum\limits_{k=1}^{4}[\partial_t^{(k)}W(x_{j+\frac{1}{2}},y_{G_y},t_n^+)]\frac{(\Delta t)^k}{k!}.
\end{align}
$W(x_{G_x},y_{j+\frac{1}{2}},t_{n+1})$ and $W(x_{i+\frac{1}{2}},y_{j+\frac{1}{2}},t_{n+1})$ can be obtained in a similar way.
  Thus the values at the next time are evaluated by the Taylor expansion (\ref{taylor2dtime}). The leading term and high order term also have been computed in the ADER method. Thus, 
  we take trivial cost to obtain the value at different times and is simpler. In fact, we can employ the Gauss-Lobatto quadrature rule to evaluate the numerical flux. Thus, the values at the next time are calculated in the ADER method and can be directly used in \eqref{derave2d} for updating the cell averages of the derivatives, which can further reduce the computational cost for our method.


%
%
%

\section{The algorithm of the ADER-SHWENO method}
\label{alg}
Assume the physical solution is given at $t=t_n$. Now the procedure of the ADER-SHWENO method is summarized:

{\bf Step 1.} Reconstruct the conservative values at the cell interface from the cell averages by the SHWENO reconstruction.

{\bf Step 2.} Express all time derivatives as functions of space derivatives by the Lax-Wendroff procedure.

{\bf Step 3.} Evaluate the values of the leading term and high order terms at the cell interfaces by solving Riemann problems.

{\bf Step 4.} Compute the numerical flux.

{\bf Step 5.} Update the conservative variables and evaluate the cell average of derivatives at the next time.

%
%
%
%
%

\section{Numerical examples}
\label{secnumexam}
\setcounter{equation}{0}
\setcounter{figure}{0}
\setcounter{table}{0}

In this section we present numerical results for a selection of one- and two-dimensional examples to demonstrate the performance of the ADER-SHWENO method proposed 
in the paper.
%
%
%
%
%
%
%
%
%
The CFL condition number is taken as 0.9 for all the computations. In addition, the classic HWENO method \cite{qiu2003} is termed as RK-HWENO in the following. And the CFL condition
number for RK-HWENO is taken as 0.6. For accuracy test, we set time step $\Delta t=O(\Delta x^{\frac{5}{3}})$ for RK-HWENO method to make sure the spatial errors dominate. 

\subsection{One-dimensional examples}

\begin{exam}{\em
\label{exam4.1}
\label{exam4.1}
We first consider Burgers' equation
\[
W_t+\left (\frac{W^2}{2}\right )_x=0,\quad x\in (0,2)
\]
subject to the initial condition $W(x,0)=0.5+\hbox{sin}(\pi x)$ and a periodic boundary condition. 

We compute the solution up to $T = \frac{0.5}{\pi}$ when the solution is still smooth and the exact solution can be computed using Newton's iteration. The errors of the ADER-SHWENO and RK-HWENO method are listed in Table \ref{ex4.1}, which shows the convergence of the fifth order. 
Moreover,  the error of the ADER-SHWENO method is smaller than ones of the RK-HWENO method.


\begin{table}
\caption{Example~\ref{exam4.1}: Solution error with periodic boundary conditions and $T=\frac{0.5}{\pi}$.}
\renewcommand{\multirowsetup}{\centering}
\begin{center}
\begin{tabular}{|c|c|c|c|c|c|c|c|c|c|c|c|c|}
\hline
method & $N$  &10      & 20 & 40 & 80 &160 &320\\
\hline
\multirow{6}{2cm}{ADER-SHWENO}
&$L^1$ &7.320e-3      & 7.184e-4 & 1.729e-5 & 6.061e-7 & 1.977e-8 & 6.198e-10\\
 &  Order      & \quad  & 3.349       & 5.377      &  4.834        & 4.938 &4.995 \\
& $L^2$   & 1.267e-2      & 1.324e-3 & 5.042e-5 &  1.738e-6  & 5.616e-8 & 1.768e-9\\
& Order   & \quad   & 3.259     & 4.714     &   4.859      & 4.951 & 4.989 \\
& $L_{\infty}$ & 3.002e-2 & 4.208e-3 & 2.637e-4  &  9.512e-6 & 3.085e-7 & 9.941e-9\\
 &Order    & \quad  & 2.835      & 3.996     &   4.793      & 4.947 & 4.955   \\
 \hline
 \multirow{6}{2cm}{RK-HWENO}

 & $L^1$  & 1.3887e-2      &1.204e-3 & 6.052e-5 &  2.294e-6 & 8.357e-8 & 2.427e-9 \\
 &  Order      & \quad & 3.525    & 4.315     & 4.722     & 4.779 & 5.106  \\
 &$L^2$  & 2.234e-2       & 2.808e-3 & 1.630e-4 &  5.974e-6 & 2.031e-7 & 6.086e-9\\
 &Order       & \quad &3.525     & 4.107      & 4.770    & 4.879 & 5.060 \\
 &$L_{\infty}$ & 9.947e-2 & 1.079e-2 & 8.033e-4  & 3.290e-5 & 1.044e-6 & 3.175e-8 \\
 &Order    & \quad       &3.204    & 3.748      & 4.610   & 4.978 &5.039   \\
 \hline
%
\end{tabular}
\end{center}
\label{ex4.1}
\end{table}

}\end{exam}

\begin{exam}{\em
\label{exam4.2}
To see the accuracy of the method for system problems, we compute the Euler equations,
\begin{equation}
\label{Euler}
\begin{pmatrix}\rho \\ \rho u \\ E \end{pmatrix}_t+\begin{pmatrix} \rho u \\ \rho u^2+P \\ u(E+P) \end{pmatrix}_x=0,
\end{equation}
where
$\rho$ is the density, $u$ is the velocity, $E$ is the energy density, and $P$ is the pressure. The equation of state is $E=\frac{P}{\gamma-1}+\frac{1}{2}\rho u^2$ with $\gamma=1.4$. The initial condition is
\[
\rho(x,0)=1+0.2\text{sin}(\pi x), \quad u(x,0)=1,\quad P(x,0)=1,
\]
and a periodic boundary condition is used. The exact solution for this problem is
$$\rho(x,t)=1+0.2\hbox{sin}(\pi(x-t)),\quad u(x,t)=1,\quad P(x,t)=1.$$

The final time is $T=10.0$. The errors of both ADER-SHWENO and RK-HWENO method in computed density 
are listed in Table~\ref{ex4.21d}. From the table one can see that the fifth order of accuracy of the scheme
is achieved for this nonlinear system. Moreover the error of ADER-SHWENO method is smaller than that with RK-HWENO method.
In addition, the accuracy will be decreased to third order for RK-HWENO if the time step is set as $\Delta t=CFL\cdot\Delta x$.


To show the efficiency of the ADER-SHWENO method, we plot the $L^{\infty}$ error as a function of the CPU time for both RK-HWENO and ADER-SHWENO methods in Fig. \ref{figtime}. The results are obtained on a single Intel Core i7 CPU with 2.5 GHz and 32.00 GB of RAM using Matlab R2021b. From the figure, one can see that the ADER-SHWENO method is more efficient than the original RK-HWENO method.

\begin{table}
\caption{Example~\ref{exam4.2}: Errors in computed density for periodic boundary conditions and $T=10.0$.}
\renewcommand{\multirowsetup}{\centering}
\begin{center}
\begin{tabular}{|c|c|c|c|c|c|c|c|c|c|c|c|c|}
\hline
method & $N$  &10      & 20 & 40 & 80 &160 \\
\hline
\multirow{6}{2cm}{ADER-SHWENO}
&$L^1$       & 1.328e-3 & 4.031e-5 & 1.274e-6 &  3.983e-8 & 1.244e-9  \\
 &  Order      & \quad    & 5.043    & 4.984     & 4.999      &  5.001        \\
& $L^2$       & 1.556e-3 & 4.545e-5 & 1.417e-6 & 4.423e-8 &  1.382e-9  \\
& Order       & \quad    & 5.097     & 5.003    & 5.002     &   5.001       \\
& $L_{\infty}$&2.263e-3  &6.976e-5  & 2.032e-6 & 6.271e-8  &  1.954e-9  \\
 &Order    & \quad       & 5.020     & 5.102     & 5.018    &   5.004        \\
 \hline
  \multirow{6}{2cm}{RK-HWENO}

 & $L^1$  & 1.276e-2     &5.794e-4 & 1.730e-5 &  5.346e-7 & 1.648e-8  \\
 &  Order      & \quad & 4.461    & 5.066    & 5.016    & 5.020   \\
 &$L^2$  & 1.375e-2      & 6.335e-4 & 1.974e-5 &  6.041e-7 & 1.862e-8 \\
 &Order       & \quad &4.440    & 5.002      & 5.033    & 5.020  \\
 &$L_{\infty}$ & 1.660e-2 & 9.258e-4 & 3.218e-5  & 1.010e-6 & 3.008e-8  \\
 &Order    & \quad       &4.164   & 4.846      & 4.994   & 5.069  \\
 \hline
\end{tabular}
\end{center}
\label{ex4.21d}
\end{table}

\begin{figure}[hbtp]
 \begin{center}
 {\includegraphics[width=15cm]{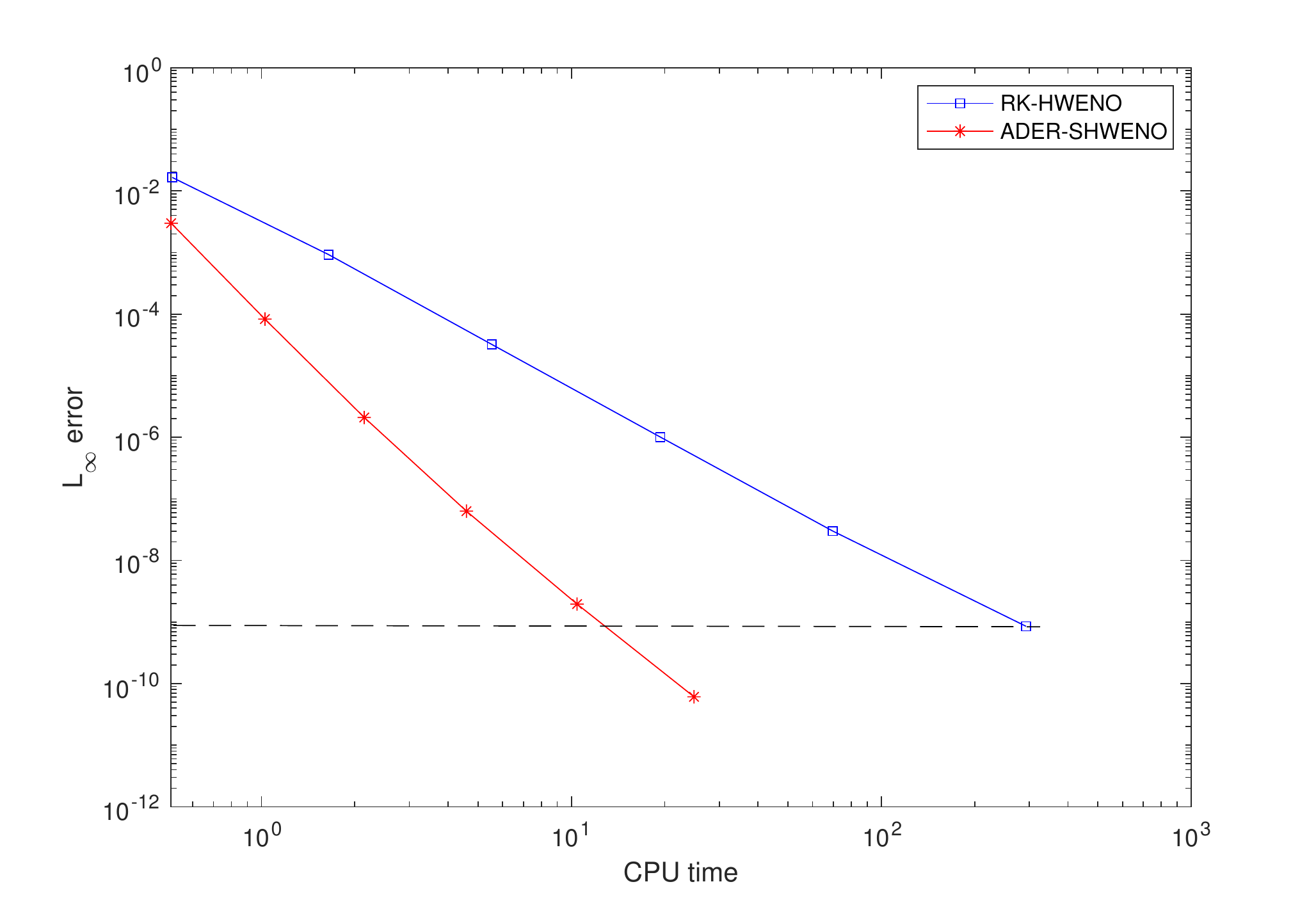}}
   \caption{$L_{\infty}$ error is plotted as a function of CPU time (in seconds).}
   \label{figtime}
   \end{center}
   \end{figure}

}\end{exam}

\begin{exam}{\em
\label{examlax}
In this example we consider the Lax problem of the Euler equations (\ref{Euler}) with the initial condition
\begin{equation}
(\rho,u,p)=
\left
\{
\begin{array}{ll}
(0.445,0.698,3.528), \quad &\text{for}\quad x<0\\
(0.5,0,0.571), \quad &\text{for}\quad x>0\notag
\end{array}
\right.
\end{equation}
and the inflow/outflow boundary condition. The computational domain is $(-5,5)$ and the integration is stopped at $T=1.3$. 

The computed density obtained by the RK-HWENO and ADER-SHWENO methods is plotted with 200 points in Fig. \ref{figlax}, in which one can see that the resolution obtained by ADER-SHWENO is comparable with that obtained by RK-HWENO.

\begin{figure}[hbtp]
 \begin{center}
 \mbox{
 {\includegraphics[width=8cm]{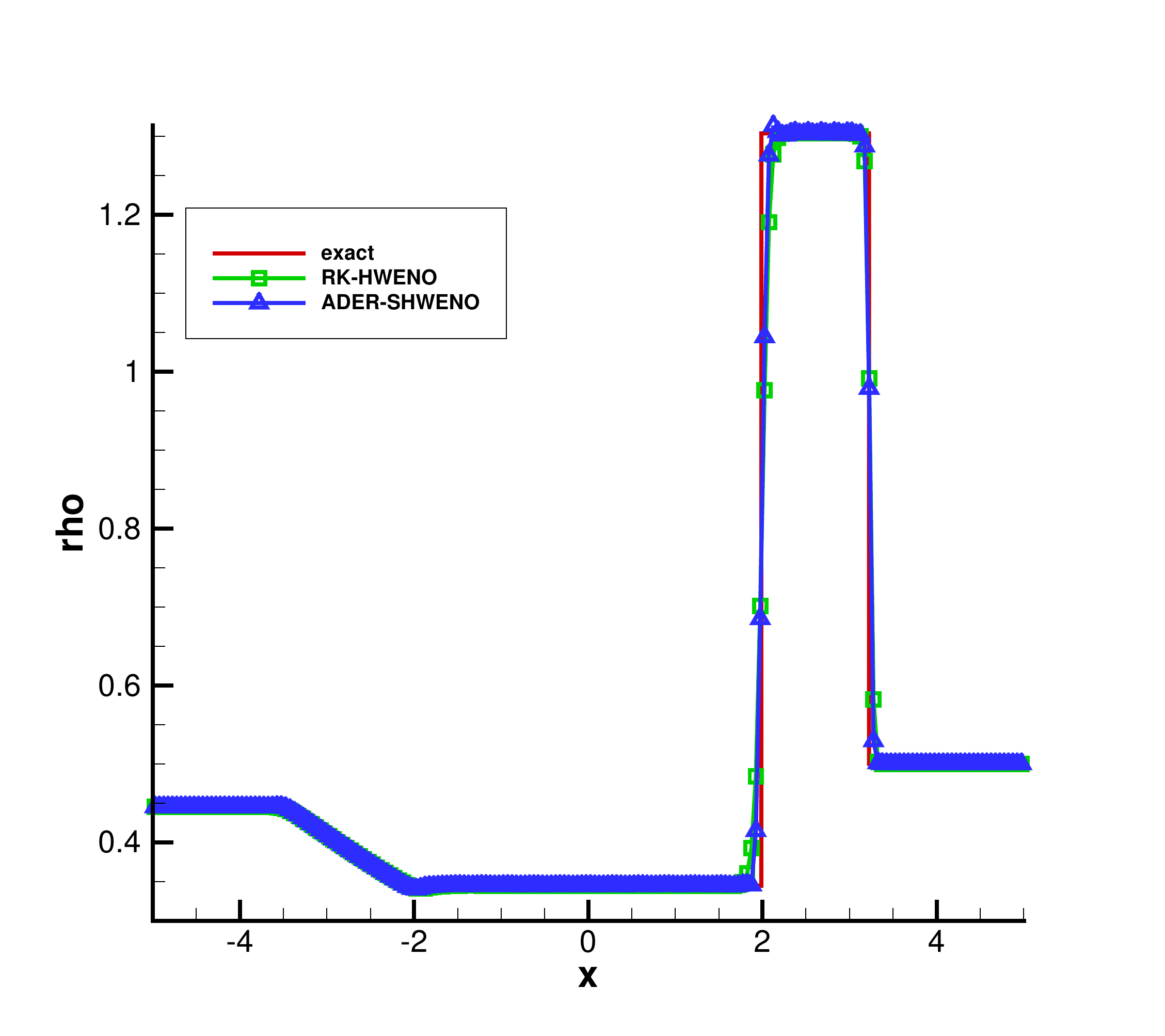}}\quad
{\includegraphics[width=8cm]{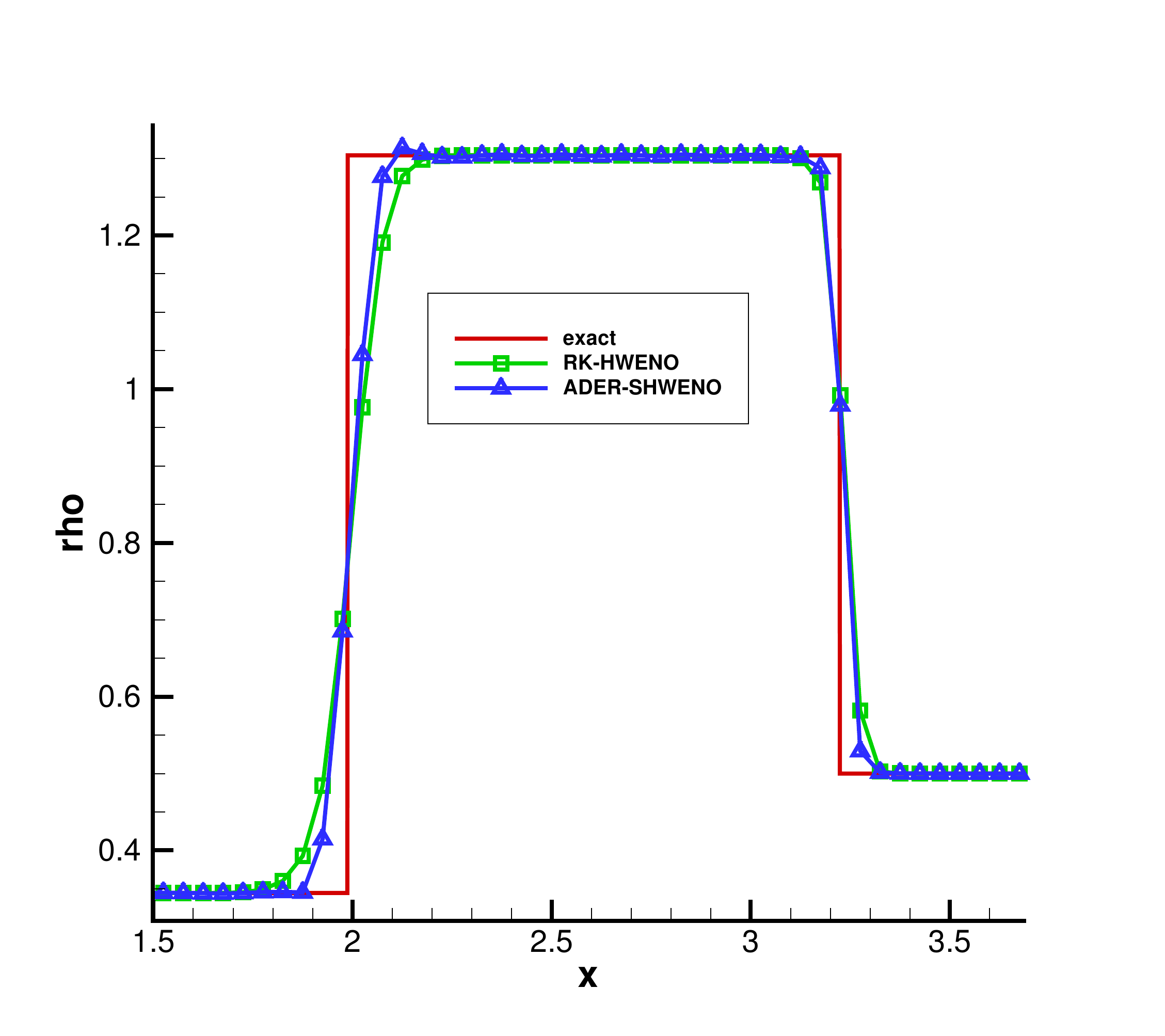}}
}
  \caption{Example~\ref{examlax} The density obtained by ADER-SHWENO method with $N=200$ uniform meshes is compared with RK-HWENO method. Green square: RK-HWENO; Blue triangle: ADER-SHWENO}
   \label{figlax}
   \end{center}
   \end{figure}

}\end{exam}

\begin{exam}{\em
\label{examshuosher}
The Shu-Osher problem is considered in this example, which contains both shocks and complex smooth region structures. We solve the Euler equations (\ref{Euler}) with a moving shock ($\hbox{Mach}=3$) interacting with  a sine wave in density. The initial condition is
\begin{equation}
(\rho,u,p)=
\begin{cases}
(3.857143,2.629369,10.333333), \quad &\text{for}\quad x<-4,\\
(1+0.2\hbox{sin}(5x),0,1), \quad &\text{for}\quad x>-4.\notag
\end{cases}
\end{equation}
The physical domain is taken as $(-5,5)$ in this computation. The computed density is shown at $T=1.8$ against an ``exact solution'' obtained by a fifth-order finite volume WENO scheme with 10,000 uniform points. 

The solution obtained by the ADER-SHWENO method with 400 uniform points is compared with the uniform mesh solutions obtained by the RK-HWENO method in Fig. \ref{figshuosher}. From the figures, one can observe that the results obtained by ADER-SHWENO method are better than ones obtained by RK-HWENO method.

\begin{figure}[hbtp]
 \begin{center}
 \mbox{
 {\includegraphics[width=8cm]{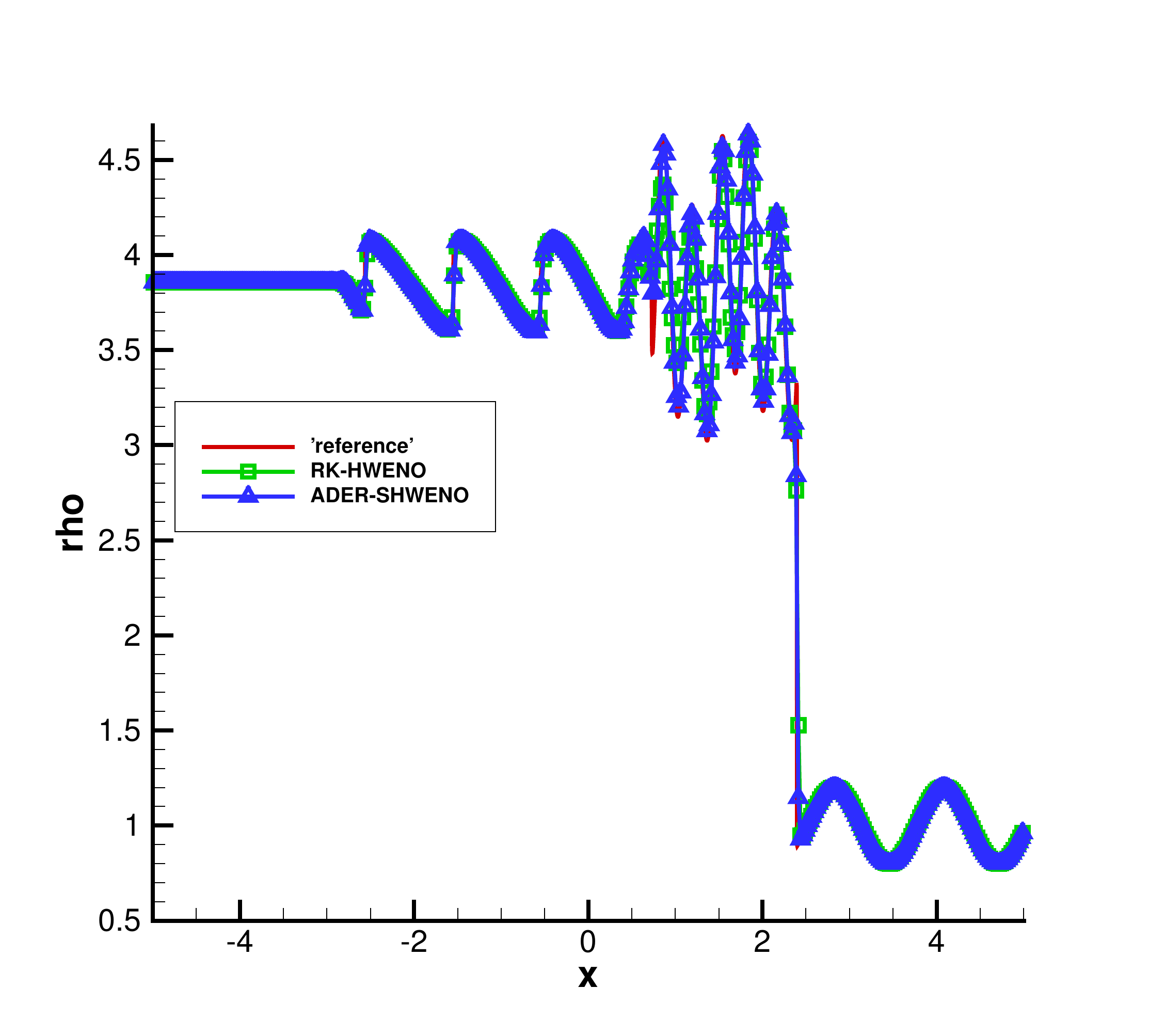}}\quad
{\includegraphics[width=8cm]{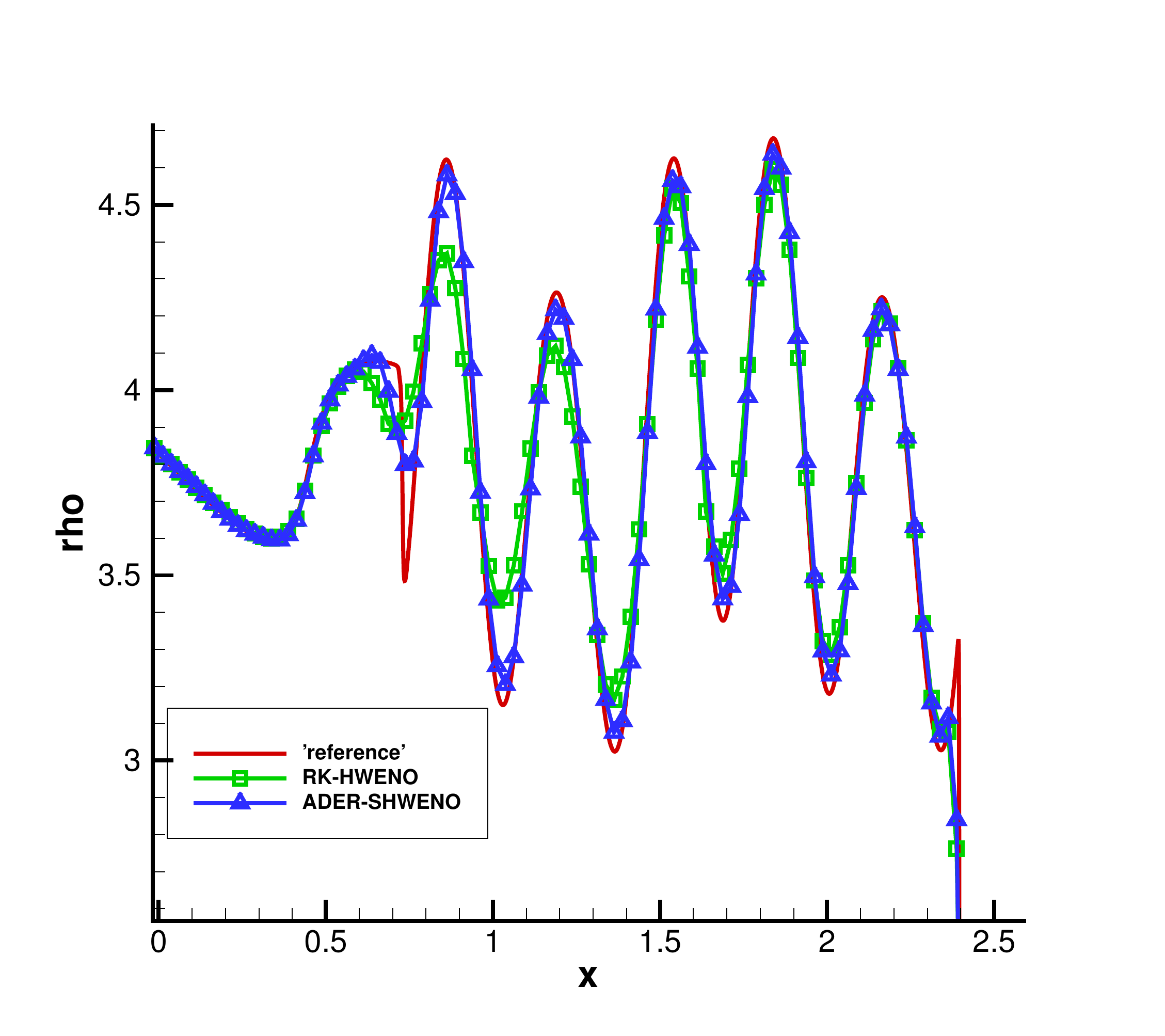}}
}

%
\caption{Example~\ref{examshuosher} The density obtained by ADER-SHWENO method with $N=400$ uniform meshes is compared with RK-HWENO method. Green square: RK-HWENO; Blue triangle: ADER-SHWENO }
   \label{figshuosher}
   \end{center}
   \end{figure}

}
\end{exam}

\begin{exam}{\em
\label{examtur}
In this example the turbulence interaction \cite{li2021} is considered. The computational domain is taken as $(-5,5)$. The initial condition is given by
\begin{equation}
(\rho,u,p)=
\begin{cases}
(1.515695,0.523346,1.80500), \quad &\text{for}\quad x<-4.5\\
(1+0.1\hbox{sin}(20\pi x),0,1), \quad &\text{for}\quad x>-4.5.\notag
\end{cases}
\end{equation}
The results of the computed density obtained by both methods with 1500 uniform meshes are plotted at $T=5$ against a ``reference solution'' obtained by the ADER-SHWENO method with 10,000 uniform points.

From Fig. \ref{figtur}, one can see that the resolution of the ADER-SHWENO method is much better than that of the RK-HWENO method, which shows the advantages of the new method for the problem containing the complex structure.

\begin{figure}[hbtp]
 \begin{center}
 \mbox{\subfigure[ADER-SHWENO]
 {\includegraphics[width=8cm]{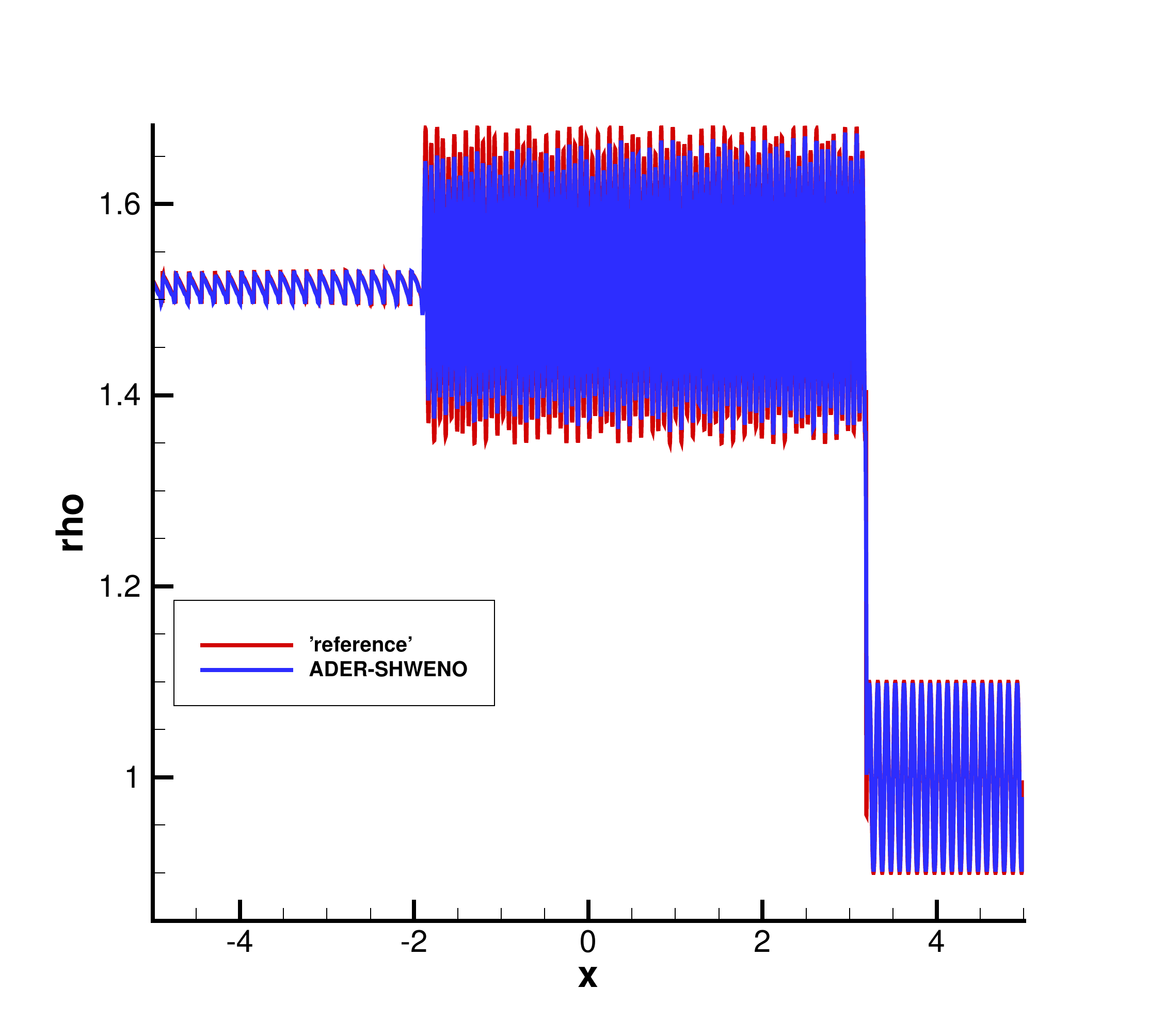}}\quad
 \subfigure[RK-HWENO]
 {\includegraphics[width=8cm]{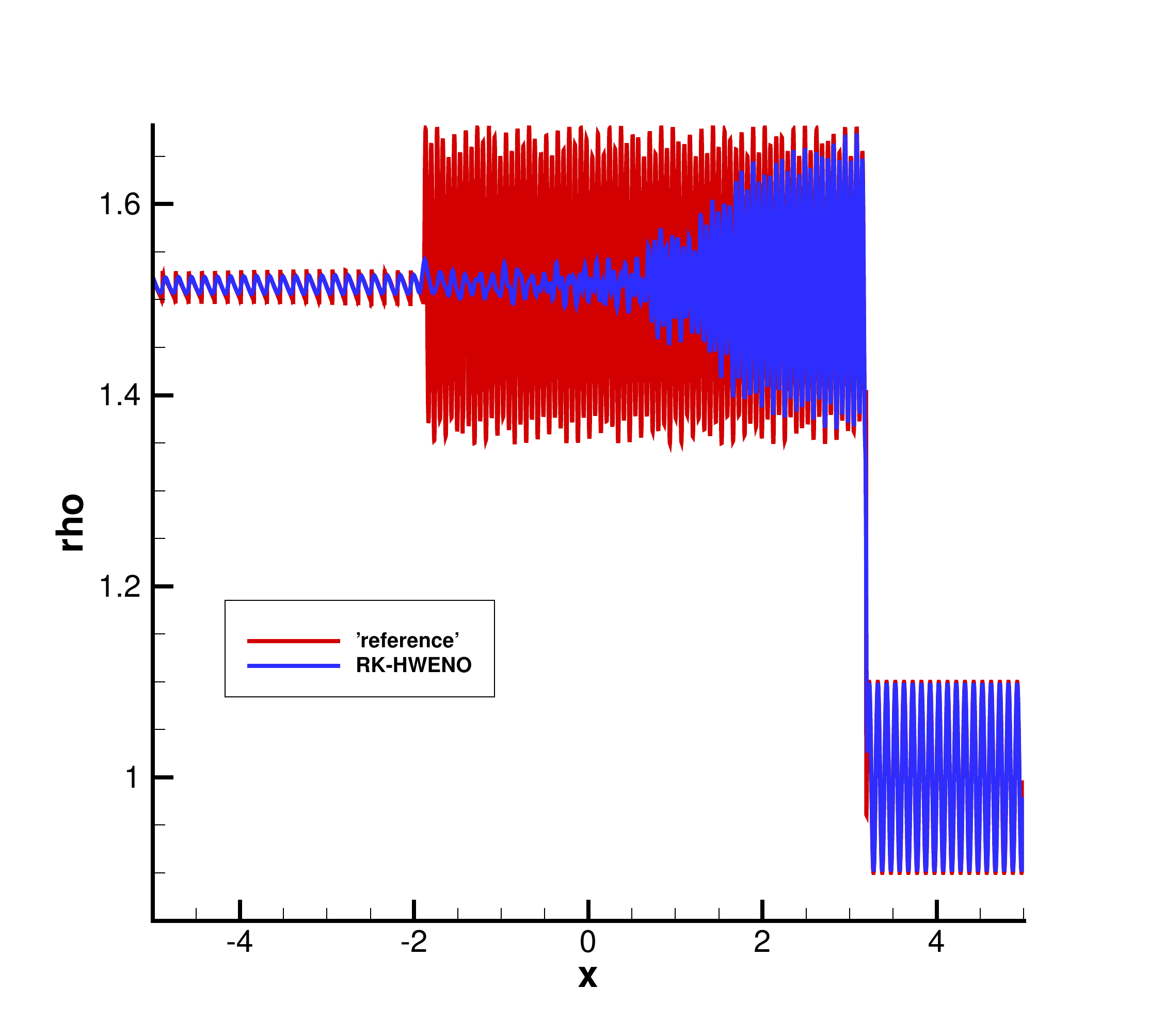}}

}
 \mbox{ \subfigure[zoom of (a)]
{\includegraphics[width=8cm]{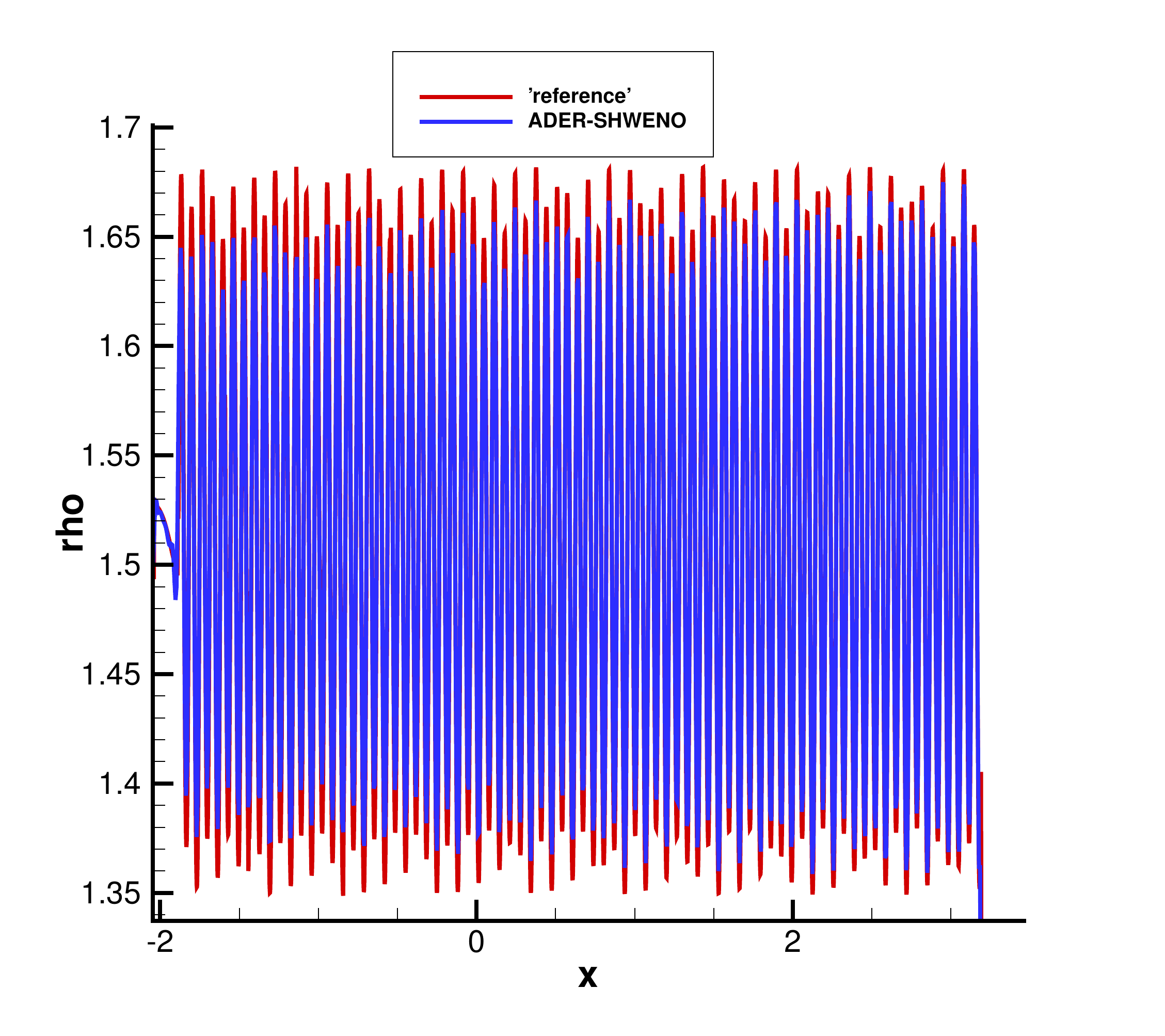}}\quad
  \subfigure[zoom of (b)]
{\includegraphics[width=8cm]{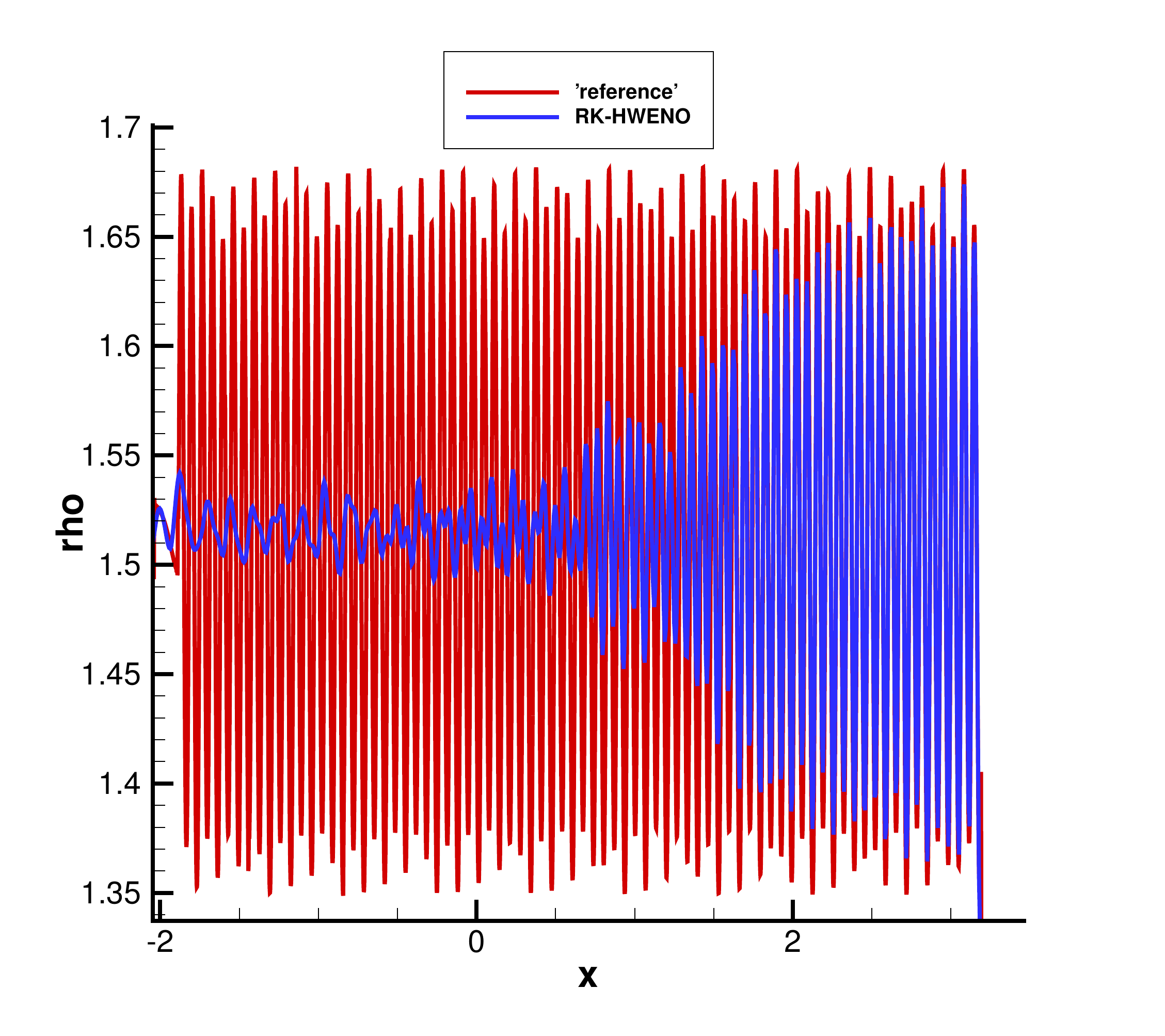}}
}

    \caption{ The density obtained by ADER-SHWENO method with $N=1500$ uniform meshes is compared with RK-HWENO method}
   \label{figtur}
   \end{center}
   \end{figure}

}
\end{exam}

\begin{exam}{\em
\label{exambw}
We consider the interaction of blast waves of the Euler equations (\ref{Euler}), which was first used by Woodward and Colella \cite{woodward1984} as a test problem for various numerical schemes.  The initial condition is given by
\begin{equation}
(\rho,u,p)=
\left
\{
\begin{array}{ll}
(1.0,0,1000), \quad& \text{for}\quad 0\leq x<0.1\\
(1.0,0,0.01), \quad& \text{for}\quad 0.1\leq x<0.9\\
(1.0,0,100),\quad&  \text{for}\quad 0.9\leq x\leq 1 .\notag
\end{array}
\right.
\end{equation}
The physical domain is taken as $(0,1)$ and a reflective boundary condition is applied to both ends. The results at time $T=0.038$ are plotted against an ``exact solution" computed by a fifth-order finite difference WENO scheme \cite{jiang1996} with 81,920 uniform mesh points in Fig. \ref{figbw}. 

From the figures we can see that the solution obtained by ADER-SHWENO method is comparable with that obtained by RK-HWENO method.

%

%
%
%
%
   
%
   \begin{figure}[hbtp]
 \begin{center}
 \mbox{
 {\includegraphics[width=8cm]{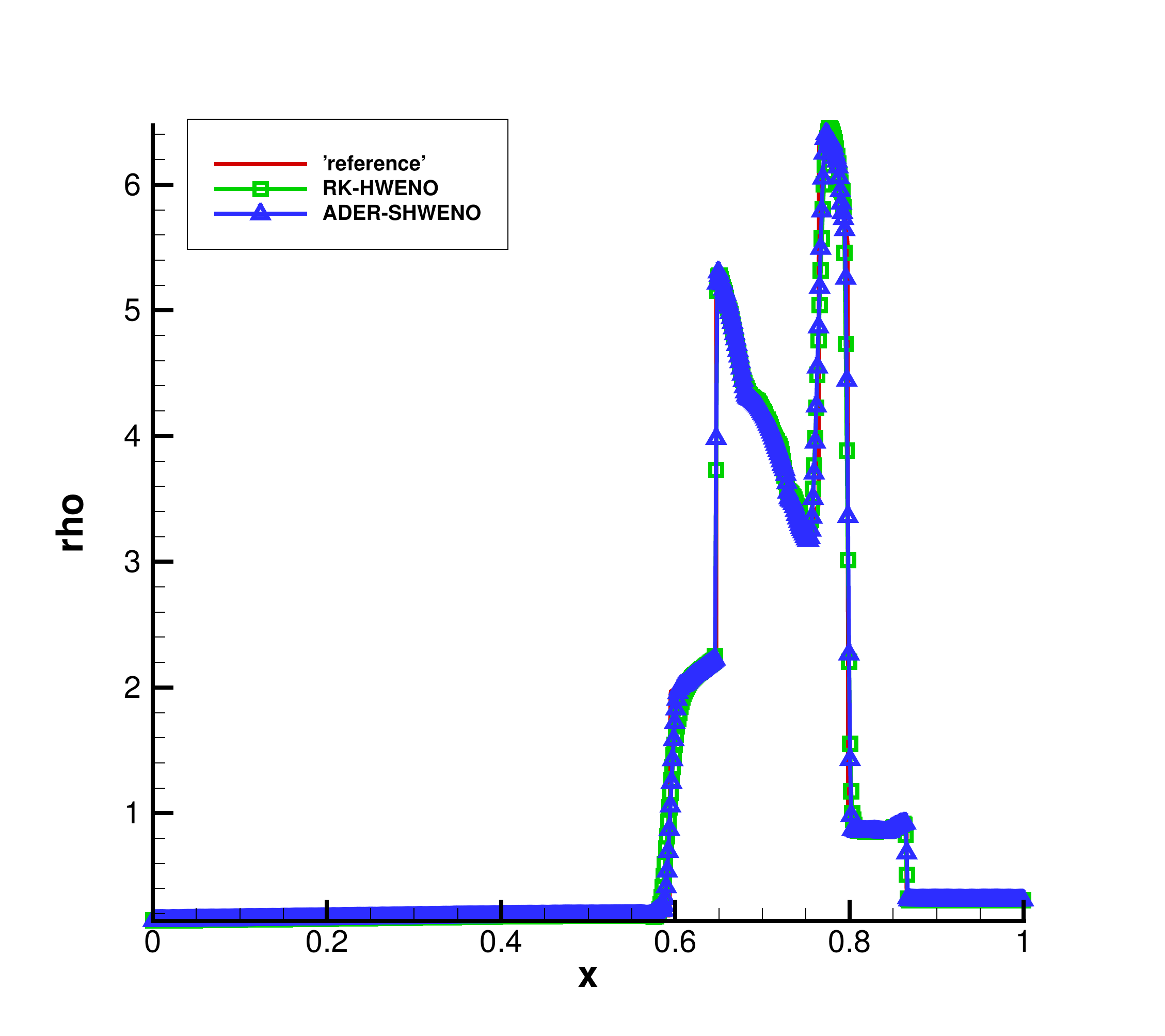}}\quad
{\includegraphics[width=8cm]{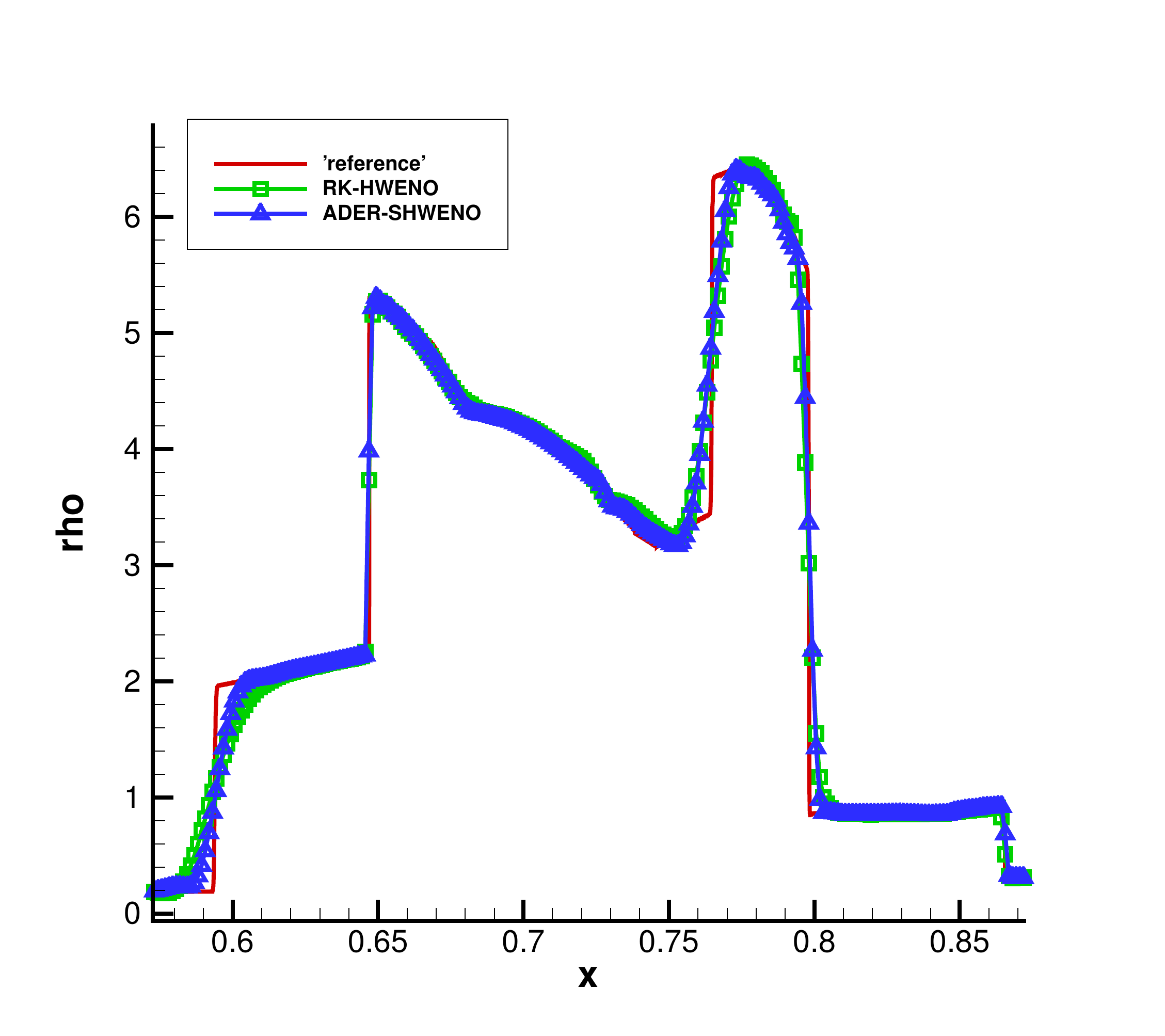}}
}

%
%
%
   \caption{Example~\ref{exambw} The density obtained by ADER-SHWENO method with $N=800$ uniform meshes is compared with RK-HWENO method. Green square: RK-HWENO; Blue triangle: ADER-SHWENO }
   \label{figbw}
   \end{center}
   \end{figure}

}
\end{exam}

\subsection{Two-dimensional examples}

%
%

\begin{exam}{\em
\label{exam4.12d}
In order to test the accuracy in two-dimensional case, we solve the Euler equations
\begin{equation}
\label{2d}
W_t+F(W)_x+G(W)_y\equiv \frac{\partial }{\partial t}
\left(
  \begin{array}{c}
    \rho  \\
     \rho\mu   \\
      \rho\nu   \\
      E \\
   \end{array}
 \right)
 +\frac{\partial }{\partial x}
 \left(
  \begin{array}{c}
    \rho\mu  \\
     \rho\mu^2+P   \\
      \rho\mu\nu   \\
      \mu(E+P) \\
   \end{array}
 \right)
 +\frac{\partial }{\partial y}
 \left(
  \begin{array}{c}
    \rho\nu  \\
     \rho\mu\nu   \\
      \rho\nu^2+P   \\
      \nu(E+P) \\
   \end{array}
 \right)
 =0,
\end{equation}
 where $\rho$ is the density, $\mu$ and $\nu$ are the velocity components
in the $x$- and $y$-direction, respectively, $E$ is the energy density, and $P$ is the pressure. The equation of state is $E=\frac{P}{\gamma-1}+\frac{1}{2}\rho(\mu^2+\nu^2)$ with $\gamma=1.4$.
 The initial condition is given by 
 $$\rho(x,y,0)=1+0.2\hbox{sin}(\pi (x+y)),\; \mu(x,y,0)=1,\; \nu(x,y,0)=1,\; P(x,y,0)=1,$$
  and a periodic boundary condition is applied in both directions.

The computational domain is $(0,2)\times (0,2)$ and the final time is $T=1$. The results in Table~\ref{ex4.12d} show the convergence of the fifth order for both ADER-SHWENO and RK-HWENO method is achieved for the Euler system in two dimensions.  In addition, the error of ADER-SHWENO  is much smaller than that of RK-HWENO. 
 

\begin{table}
\caption{Example~\ref{exam4.12d}: Solution error with periodic boundary conditions and $T=1$.}
\begin{center}
\begin{tabular}{|c|c|c|c|c|c|c|c|c|c|c|c|c|}
\hline
method & $N\times M$  &$10\times 10$      & $20\times20$ & $40\times40$ & $80\times80$ &$160\times160$ \\
\hline
\multirow{6}{2cm}{ADER-SHWENO}
   &$L^1$ &8.405e-4      & 8.659e-6 & 2.584e-7 & 7.873e-9  & 2.427e-10\\
  &  Order      & \quad  & 6.601      & 5.067      &  5.037        & 5.020 \\
 & $L^2$   & 1.101e-3      & 9.391e-6 & 2.833e-7 &  8.692e-9  & 2.689e-10 \\
 & Order   & \quad   & 6.873     & 5.051      &   5.026      & 5.015 \\
 & $L_{\infty}$ & 1.887e-3 & 1.558e-5 & 4.080e-7  &  1.232e-8 & 3.820e-10 \\
 & Order    & \quad  & 6.920      & 5.255     &   5.049      & 5.011   \\
 \hline
 \multirow{6}{2cm}{RK-HWENO}
 &  $L^1$ &9.860e-3      & 3.836e-4 & 1.106e-5 & 3.344e-7  & 1.024e-8\\
  &  Order      & \quad  & 4.684      & 5.116     &  5.048        & 5.029 \\
  &$L^2$   & 1.036e-2      & 4.078e-4 & 1.231e-5 &  3.735e-7  & 1.143e-8 \\
 & Order   & \quad   & 4.667     & 5.050      &   5.043      & 5.030 \\
  &$L_{\infty}$ &1.340e-2 & 5.734e-4 & 1.983e-5  &  6.213e-7 & 1.820e-8 \\
  &Order    & \quad  & 4.547      & 4.854     &   4.996      & 5.093   \\
 \hline
 \end{tabular}
\end{center}
\label{ex4.12d}
\end{table}


}\end{exam}

\begin{exam}{\em
\label{isentro}
In this example, the two-dimensional isotropic vortex problem is considered. The computational domain is taken as $(0,10)\times (0,10)$ and the initial condition is 
\begin{equation*}
\rho= (1-\frac{25(\gamma -1)}{8\gamma \pi^2}e^{1-r^2})^{\frac{1}{\gamma-1}},\\
\mu=1-\frac{5}{2\pi}e^{\frac{1-r^2}{2}}(y-5),\\
\nu=1+\frac{5}{2\pi}e^{\frac{1-r^2}{2}}(x-5),\\
P=\rho^{\gamma},
\end{equation*}
where $r^2=(x-5)^2+(y-5)^2$. The periodic boundary condition is employed in both directions. The exact solution is the vortex along the upper right direction with velocities $(\mu,\nu)=(1,1)$. We compute the numerical solution at the output time $T=10$. At this time the vortex returns to the initial position. 

The errors of the computed density are listed in Table~\ref{isotro} for both methods, in which one can observe that both ADER-SHWENO method 
can achieve the fifth order accuracy for this nonlinear smooth problem. 
From the table one can see that the error of ADER-SHWENO is much smaller the same as the previous example.

%
\begin{table}
\caption{Example~\ref{isentro}: Solution error with periodic boundary conditions and $T=10$.}
\begin{center}
\begin{tabular}{|c|c|c|c|c|c|c|c|c|c|c|c|c|}
\hline
method&$N\times M$  &$10\times 10$      & $20\times20$ & $40\times40$ & $80\times80$ &$160\times160$  &$320\times320$ \\
\hline
\multirow{6}{2cm}{ADER-SHWENO}
 & $L^1$ &1.870e-2      & 2.834e-3 & 1.113e-4 & 4.073e-6  & 1.329e-7 & 4.693e-9 \\
 &  Order      & \quad  & 2.722     & 4.670      &  4.771        & 4.938  & 4.824\\
& $L^2$   & 4.286e-2     & 5.972e-3 & 2.879e-4 &  1.160e-5  & 3.835e-7 & 1.219e-8 \\
& Order   & \quad   & 2.843     & 4.375     &   4.633      & 4.919             & 4.976 \\
& $L_{\infty}$ & 3.447e-1 & 5.292e-2 & 1.883e-3  &  7.588e-5 & 2.506e-6 & 7.822e-8\\
& Order    & \quad  & 2.703      & 4.813    &   4.633      & 4.920    & 5.002 \\
 \hline
 \multirow{6}{2cm}{RK-HWENO}
 & $L^1$ &2.609e-2      & 1.442e-2 & 7.669e-3 & 4.481e-3  & 2.074e-4 & {1.123e-6} \\
 &  Order      & \quad  & 0.855     & 0.911            &  0.775       & 4.433  & 7.529\\
& $L^2$   & 5.638e-2     & 3.372e-2 & 1.666e-2   &  1.057e-2  & 4.901e-4 & 1.719e-6 \\
& Order   & \quad             & 0.742     & 1.017       &   0.656     & 4.431        & 8.155 \\
& $L_{\infty}$ & 4.389e-1 & 2.506e-1 & 1.348e-1  &  8.709e-2 & 8.295e-3 & 1.171e-5\\
& Order    & \quad  & 0.809      & 0.895    &   0.630     & 3.392    & 9.468 \\
 \hline
 \end{tabular}
\end{center}
\label{isotro}
\end{table}


}
\end{exam}

\begin{exam}{\em
\label{examrp}
To show the performance of our method for Riemann problem genuinely in two dimensions, two 2D Riemann problems are considered in the example. The first one we consider is a 2D Riemann problem with shock waves \cite{lax1998} with the initial condition 
\begin{equation*}
(\rho,\mu,\nu,P)=
\begin{cases}
(1.1,0.0,0.0,1.1),  \; & \text{if}\;x > 0.5,\; y>0.5,\\
(0.5065,0.8939,0.0,0.35), \;  &\text{if}\;x < 0.5,\; y>0.5,\\
(1.1,0.8939,0.8939,1.1), \; &\text{if}\;x <0.5,\; y<0.5,\\
(0.5065,0.0,0.8939,0.35), \; &\text{if}\;x > 0.5,\; y<0.5,\\
\end{cases}
\end{equation*}
which is the case of left forward shock, right backward shock, upper backward shock and lower forward shock in \cite{lax1998}.
And the non-reflecting boundary conditions are employed on all the boundaries. The computational domain is taken as $(0,1)\times (0,1)$ and the stop time is taken as $0.25$. 

The computed density contours using both methods with $200\times 200$ and $400\times 400$ uniform meshes are plotted in Fig. \ref{figshock}, which shows that the results obtained by the refiner mesh have a better resolution. Besides,  one can observe that the results are comparable for both methods from the figures,.

\begin{figure}[hbtp]
 \begin{center}
 \mbox{\subfigure[ADER-SHWENO, $N=200\times 200$]
 {\includegraphics[width=8cm]{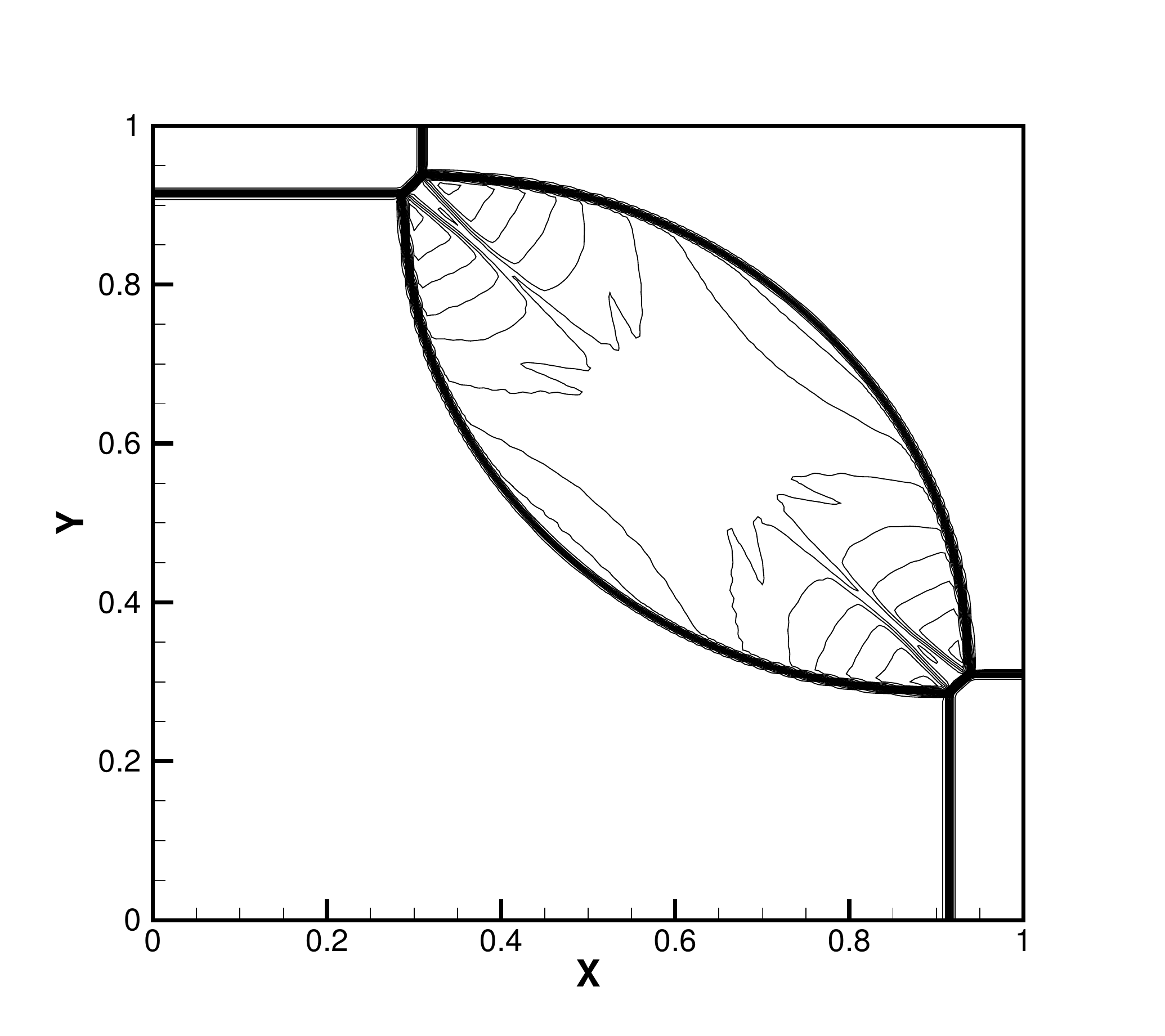}}\quad
 \subfigure[RK-HWENO, $N=200\times 200$]
 {\includegraphics[width=8cm]{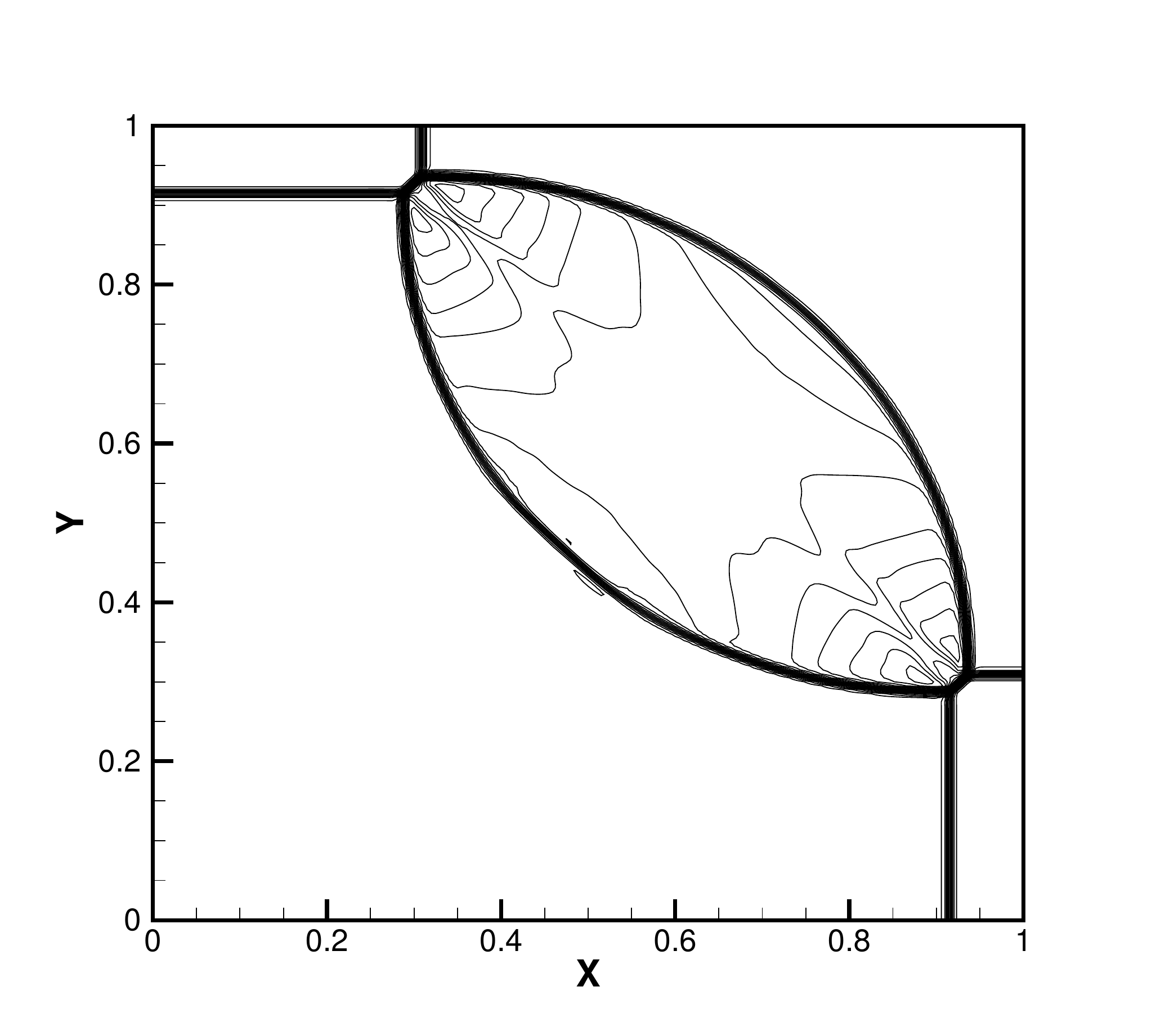}}

   }
   \mbox{ \subfigure[ADER-SHWENO, $N=400\times 400$]
{\includegraphics[width=8cm]{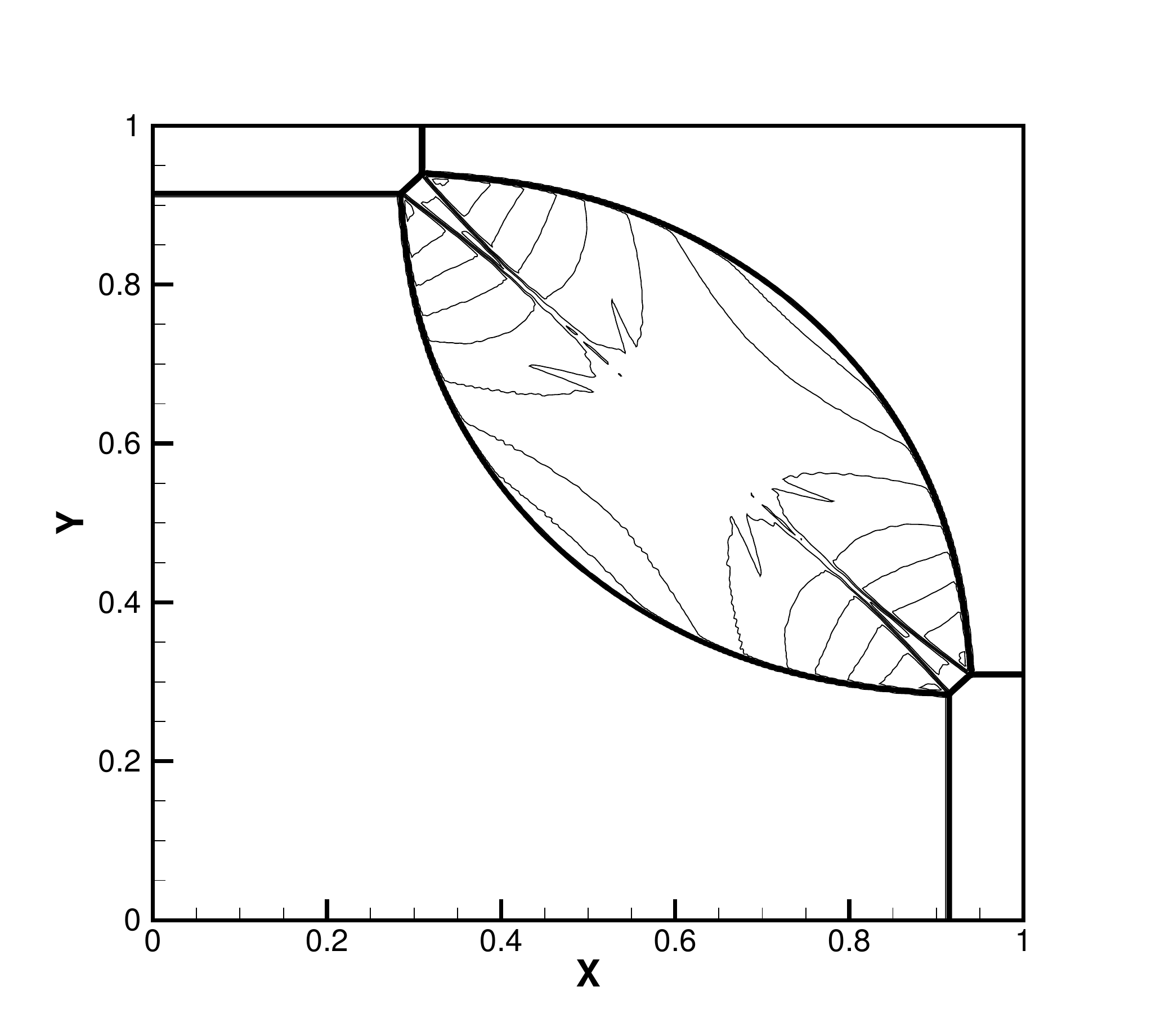}}\quad
  \subfigure[RK-HWENO, $N=400\times 400$]
{\includegraphics[width=8cm]{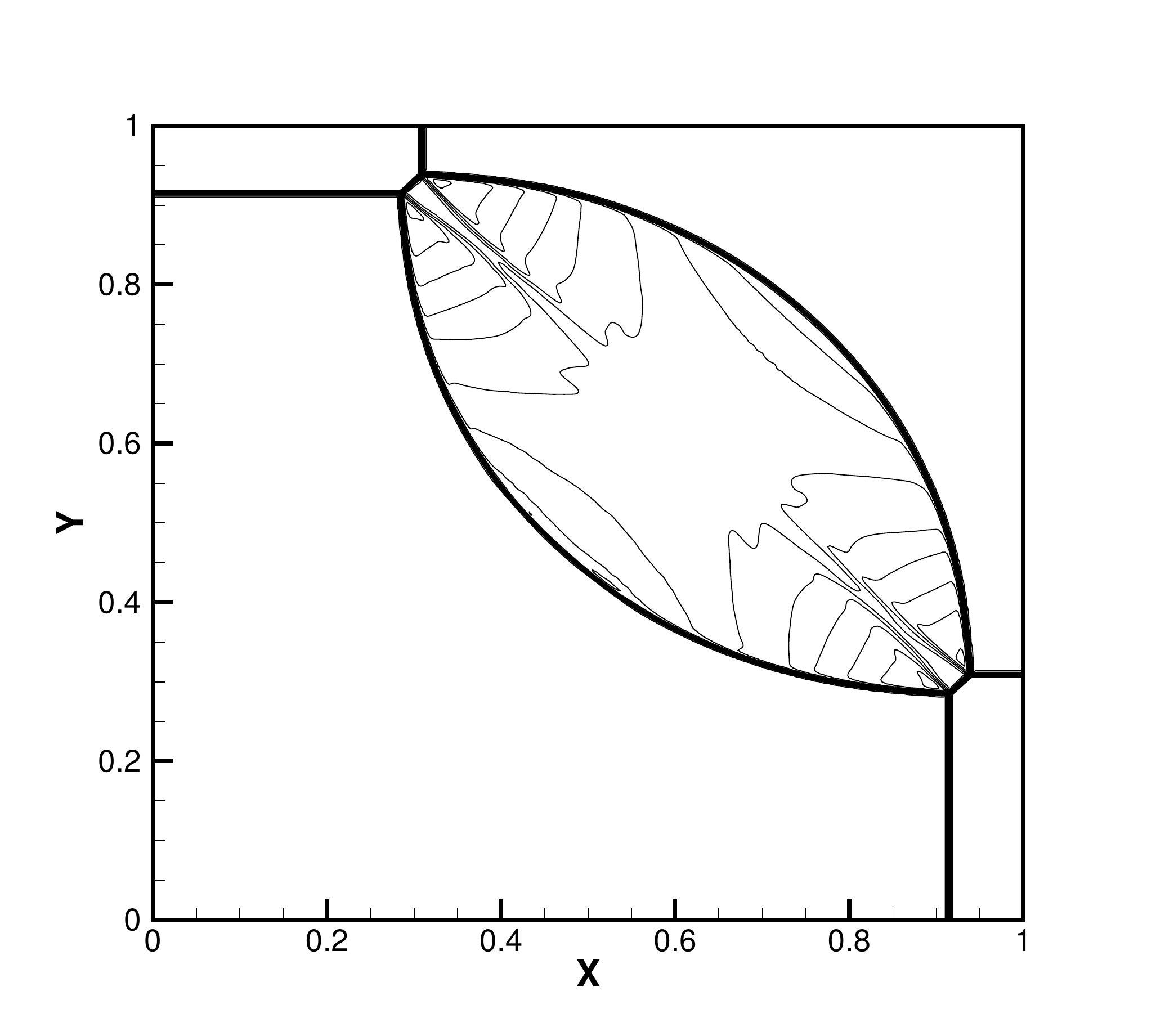}}
   }

   \caption{Example~\ref{examrp} 30 density contours from $0.5$ to 1.9.}
   \label{figshock}
   \end{center}
   \end{figure}
   
   The second problem is the shear instabilities among four initial contact discontinuities \cite{lax1998,li2021}. The initial condition is given by
\begin{equation*}
(\rho,\mu,\nu,P)=
\begin{cases}
(1.0,0.75,-0.5,1.0),  \; & \text{if}\;x > 0.5,\; y>0.5,\\
(2.0,0.5 ,  0.5,1.0), \;  &\text{if}\;x < 0.5,\; y>0.5,\\
(1.0,-0.75,0.5,1.0), \; &\text{if}\;x <0.5,\; y<0.5,\\
(3.0,-0.75,-0.5,1.0), \; &\text{if}\;x > 0.5,\; y<0.5.\\
\end{cases}
\end{equation*}
And the non-reflecting boundary conditions are adopted on all the boundaries. The computational domain is taken as $(0,1)\times (0,1)$. 

Fig. \ref{figcontact} is the result of the density computed by the ADER-SHWENO  and RK-HWENO methods with $200\times 200$ and $400\times 400$  uniform meshes at $T=0.8$. From the figure, one can observe that more complex  structures are captured for the refiner mesh. Moreover, the ADER--SHWENO method captures more complex structures than the RK-HWENO method, which indicates the advantages of the newly-developed method.

\begin{figure}[hbtp]
 \begin{center}
 \mbox{\subfigure[ADER-SHWENO, $N=200\times 200$]
 {\includegraphics[width=7cm]{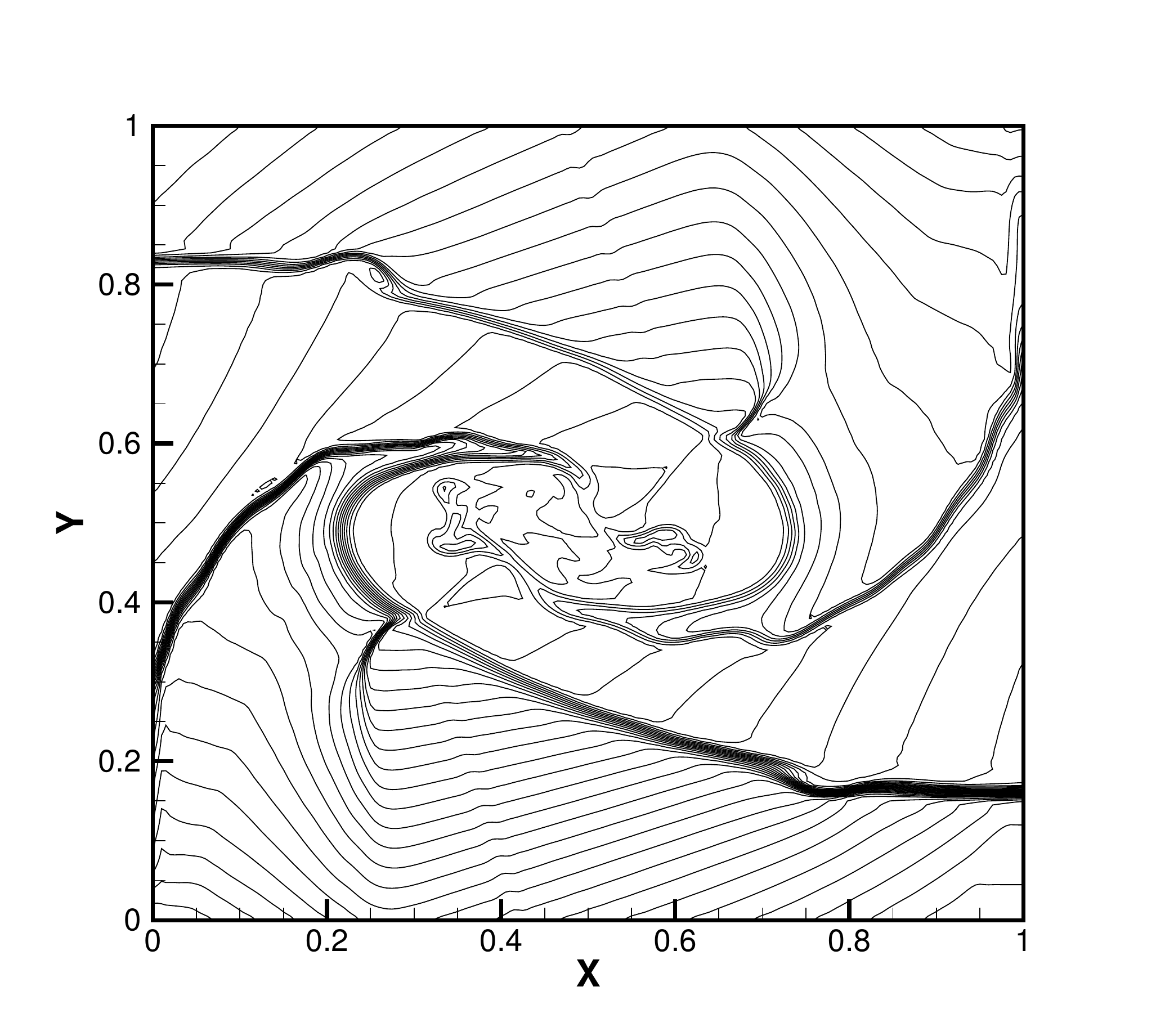}}\quad
 \subfigure[RK-HWENO, $N=200\times 200$]
 {\includegraphics[width=7cm]{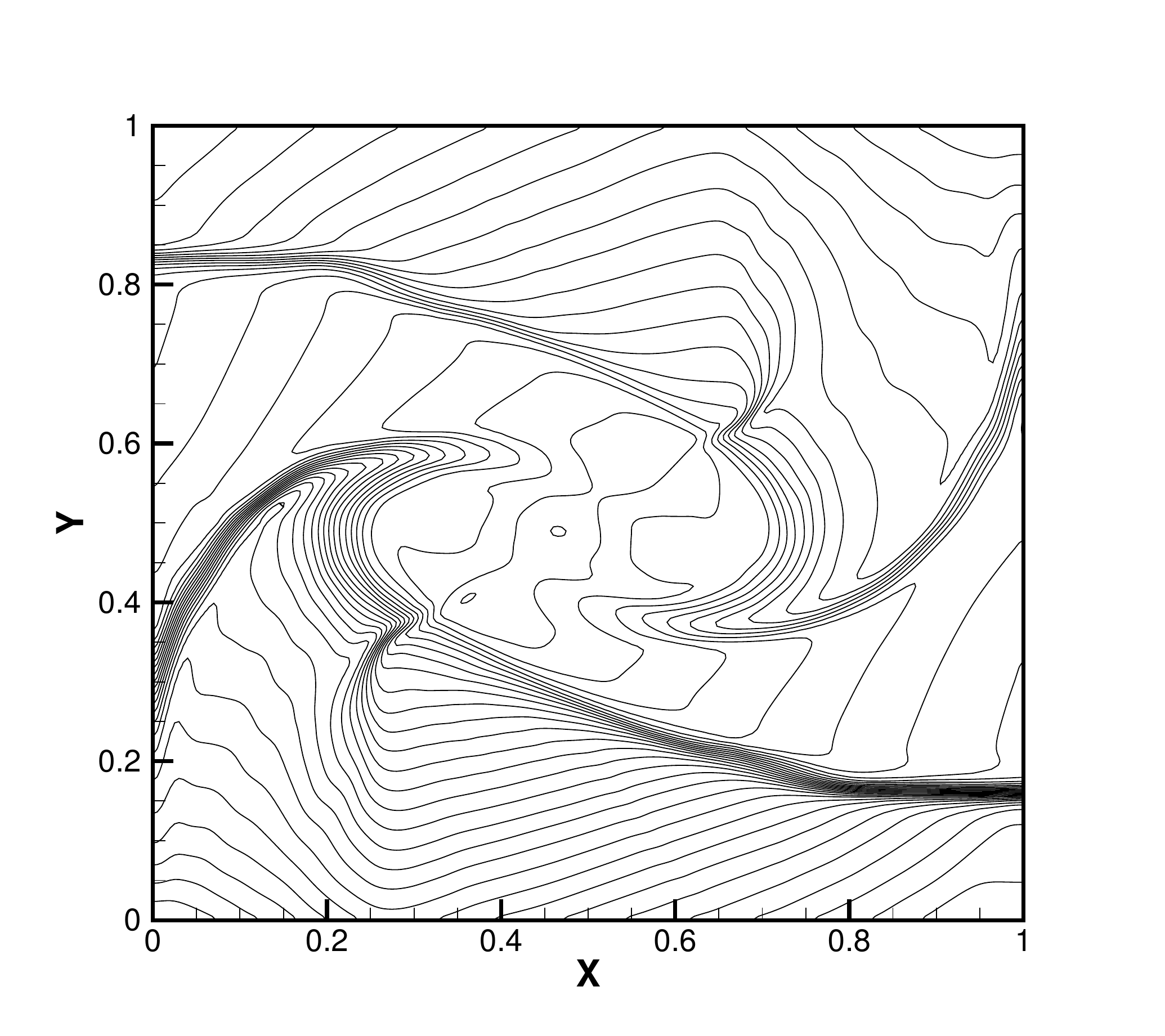}}
   }
   \mbox{
     \subfigure[ADER-SHWENO, $N=400\times 400$]
{\includegraphics[width=7cm]{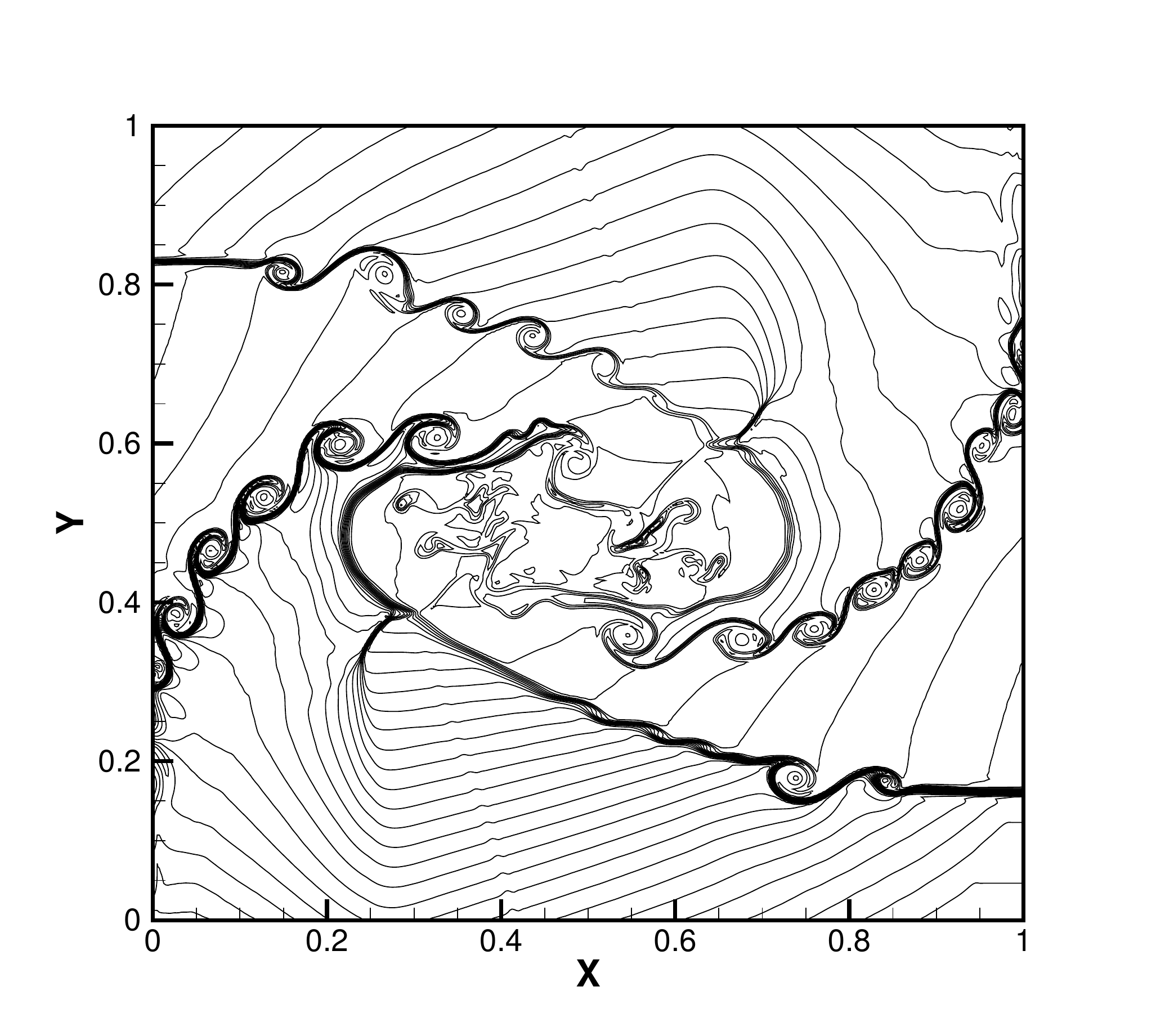}}\quad
  \subfigure[RK-HWENO, $N=400\times 400$]
{\includegraphics[width=7cm]{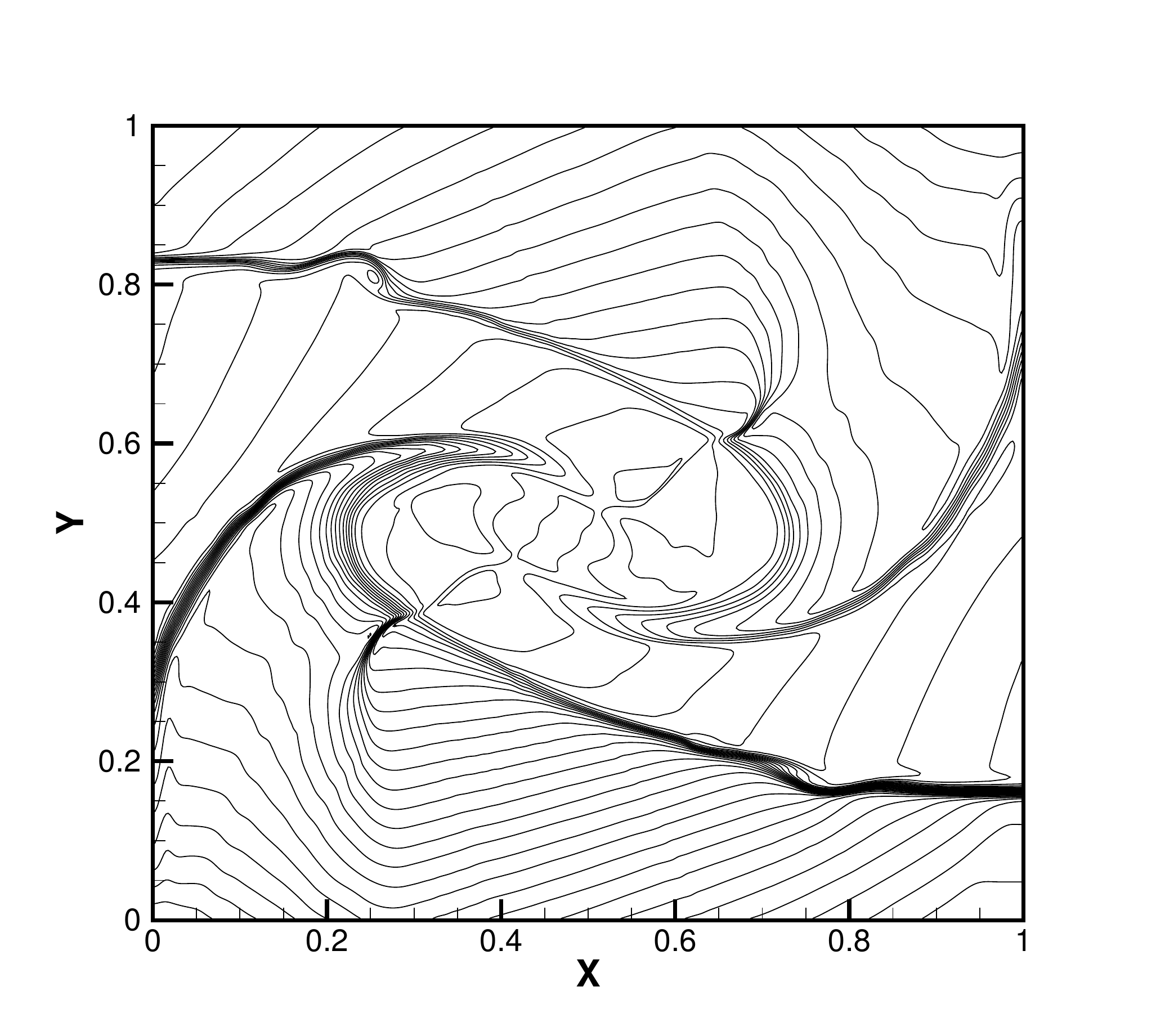}}
   }
   \caption{Example~\ref{examrp} 30 density contours from $0.2$ to 2.4.}
   \label{figcontact}
   \end{center}
   \end{figure}

}
\end{exam}

%
%
%
%

\begin{exam}{\em
\label{doublemach}

This is the double Mach reflection problem \cite{woodward1984}. We solve the Euler equations (\ref{2d}) in a computational domain of $(0,4)\times (0,1)$. The initial condition is given by 
\begin{equation}
W=
\left\{
\begin{array}{ll}
(8,57.1597,-33.0012,563.544)^T, \quad &\text{for} \quad  y\geq h(x,0)\\
(1.4,0,0,2.5)^T, \quad &\text{otherwise}
\end{array}
\right.
\notag
\end{equation}
where $h(x,t)=\sqrt{3}(x-\frac{1}{6})-20t$.
The exact post shock condition is imposed from $0$ to $\frac{1}{6}$ at the bottom while the reflection boundary condition for the rest of the bottom boundary. At the top, the boundary condition is the values that describe the exact motion of the Mach $10$ shock. On the left and right boundaries, the inflow and outflow boundary conditions are used, respectively. The final time is $T=0.2$. 

The density contours are shown in Fig. \ref{figdm} on $(0,3)\times (0,1)$.  The complex regions are plotted in Fig. \ref{figdmzoom}. From the figures, one can see that the resolution is improved as the meshes are refined. In addition, more complex vortex are captured by the ADER-SHWENO method.

\begin{figure}[hbtp]
 \begin{center}
 \mbox{\subfigure[ADER-SHWENO, $N=480\times 120$]
 {\includegraphics[width=8cm]{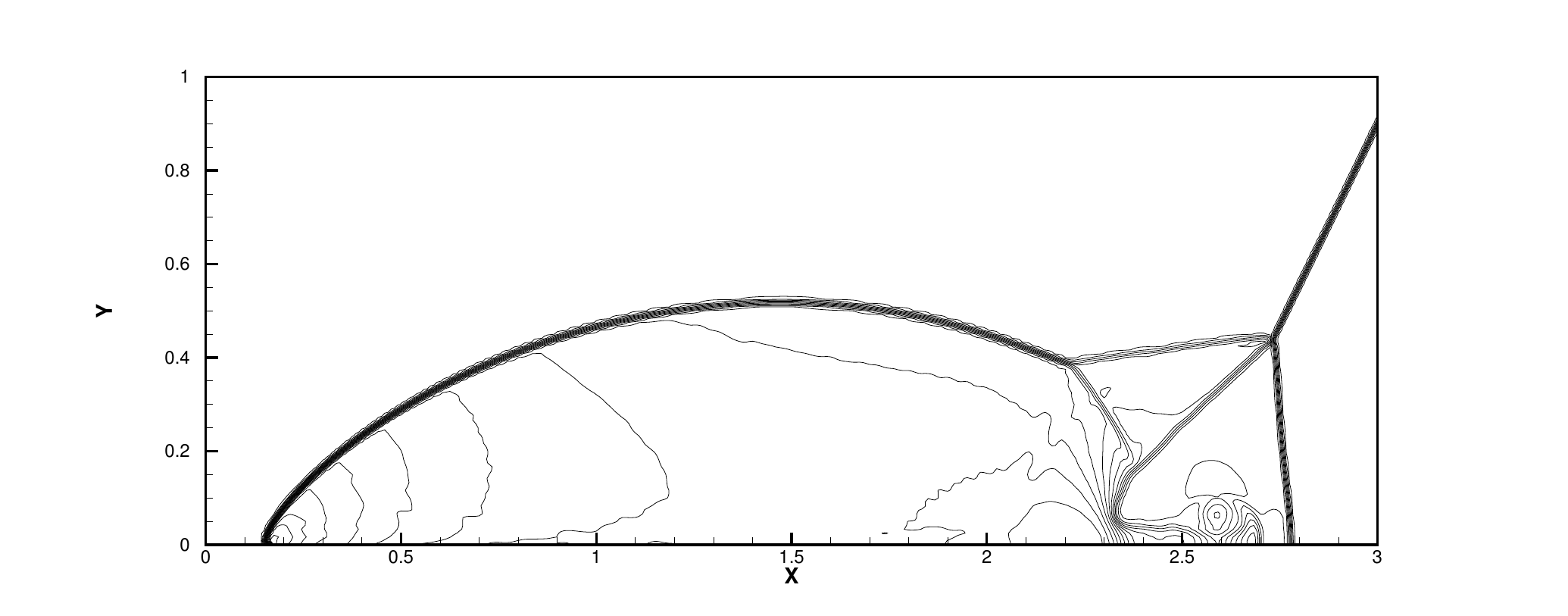}}\quad
 \subfigure[RK-HWENO, $N=480\times 120$]
 {\includegraphics[width=8cm]{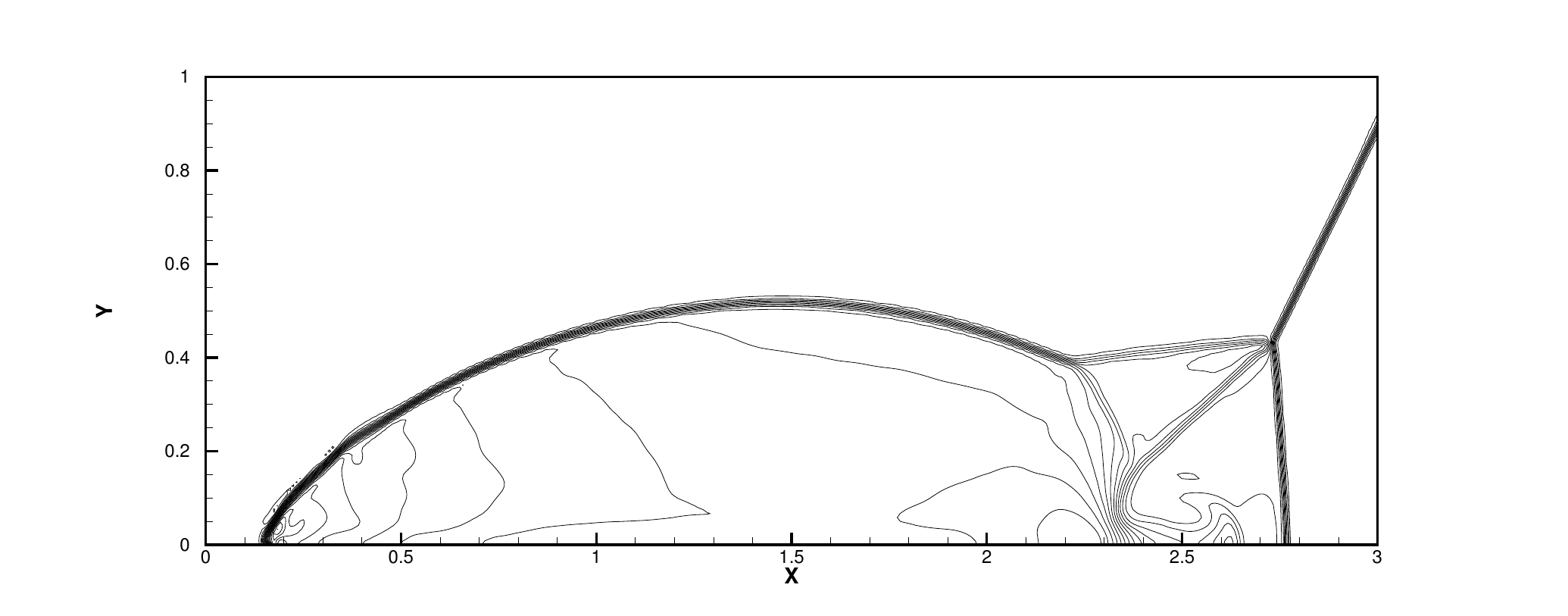}}
}
 \mbox{ \subfigure[ADER-SHWENO,  $N=960\times 240$]
{\includegraphics[width=8cm]{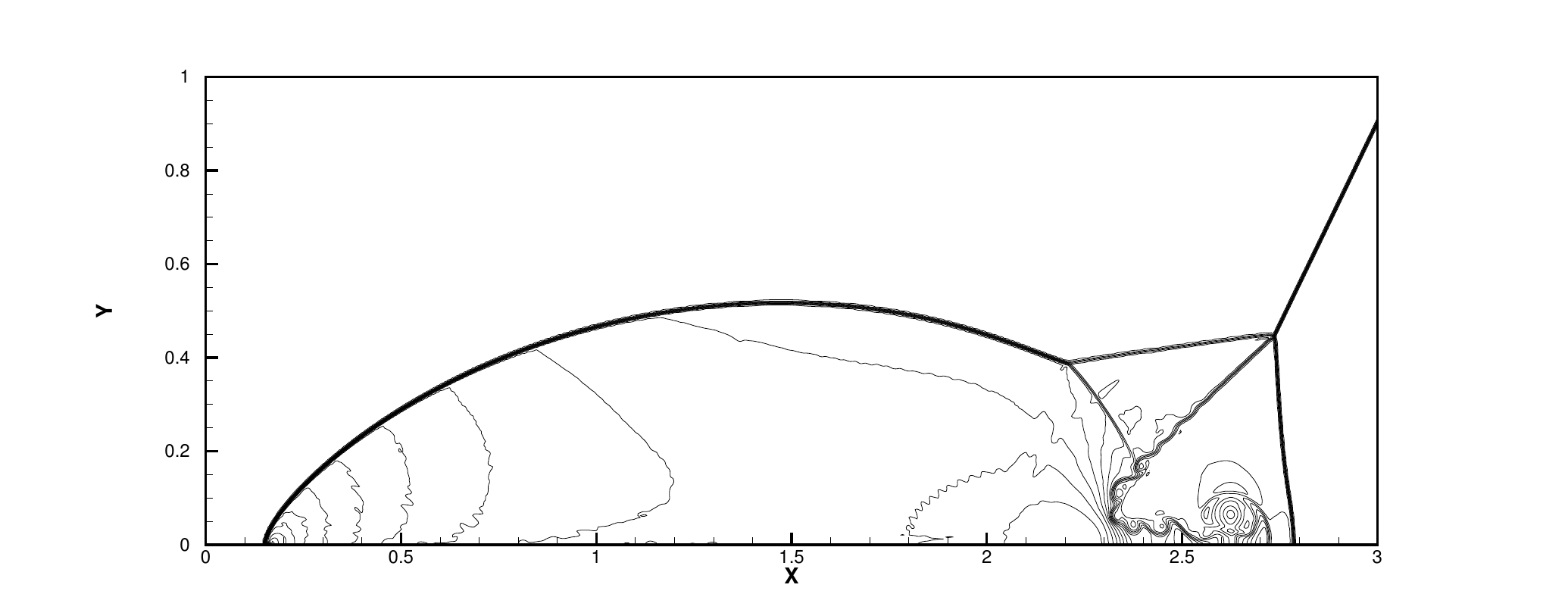}}\quad
 \subfigure[RK-HWENO,  $N=960\times 240$]
{\includegraphics[width=8cm]{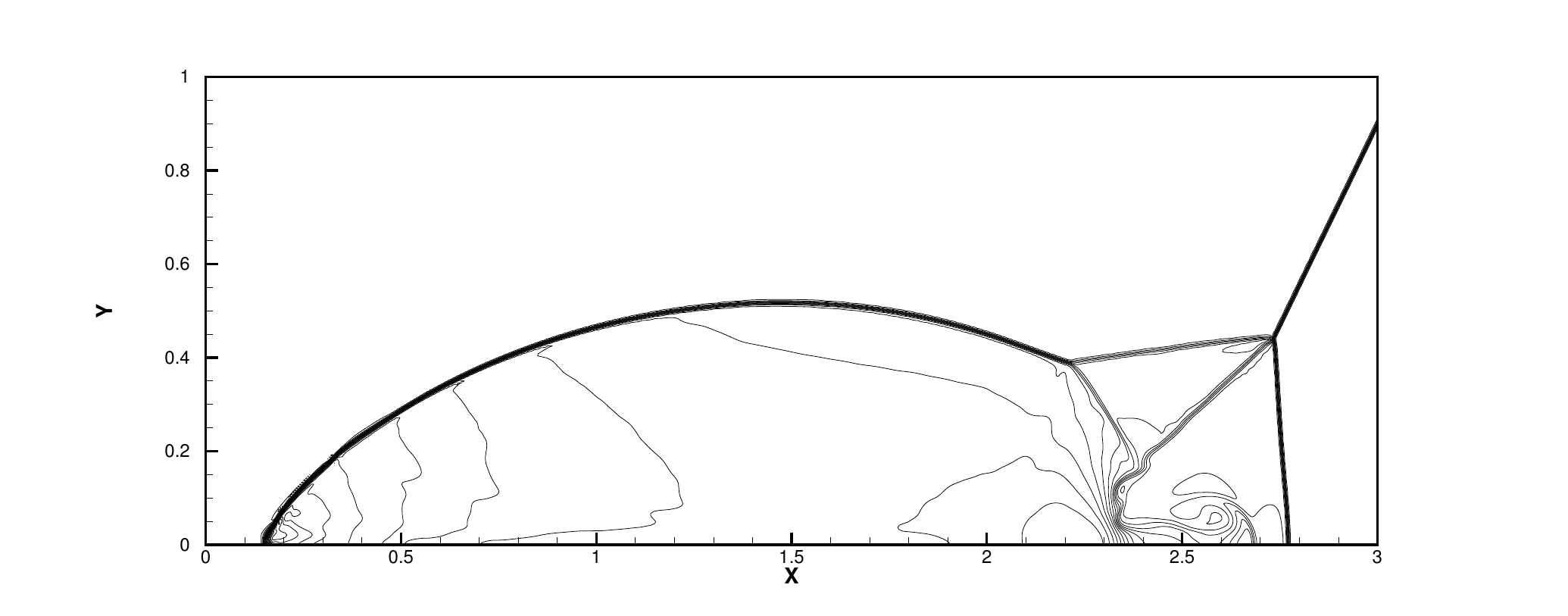}}
}


   \caption{Example~\ref{doublemach} 30 density contours from $2$ to 22.}
   \label{figdm}
   \end{center}
   \end{figure}
   
   \begin{figure}[hbtp]
 \begin{center}
 \mbox{\subfigure[ADER-SHWENO, $N=480\times 120$]
 {\includegraphics[width=8cm]{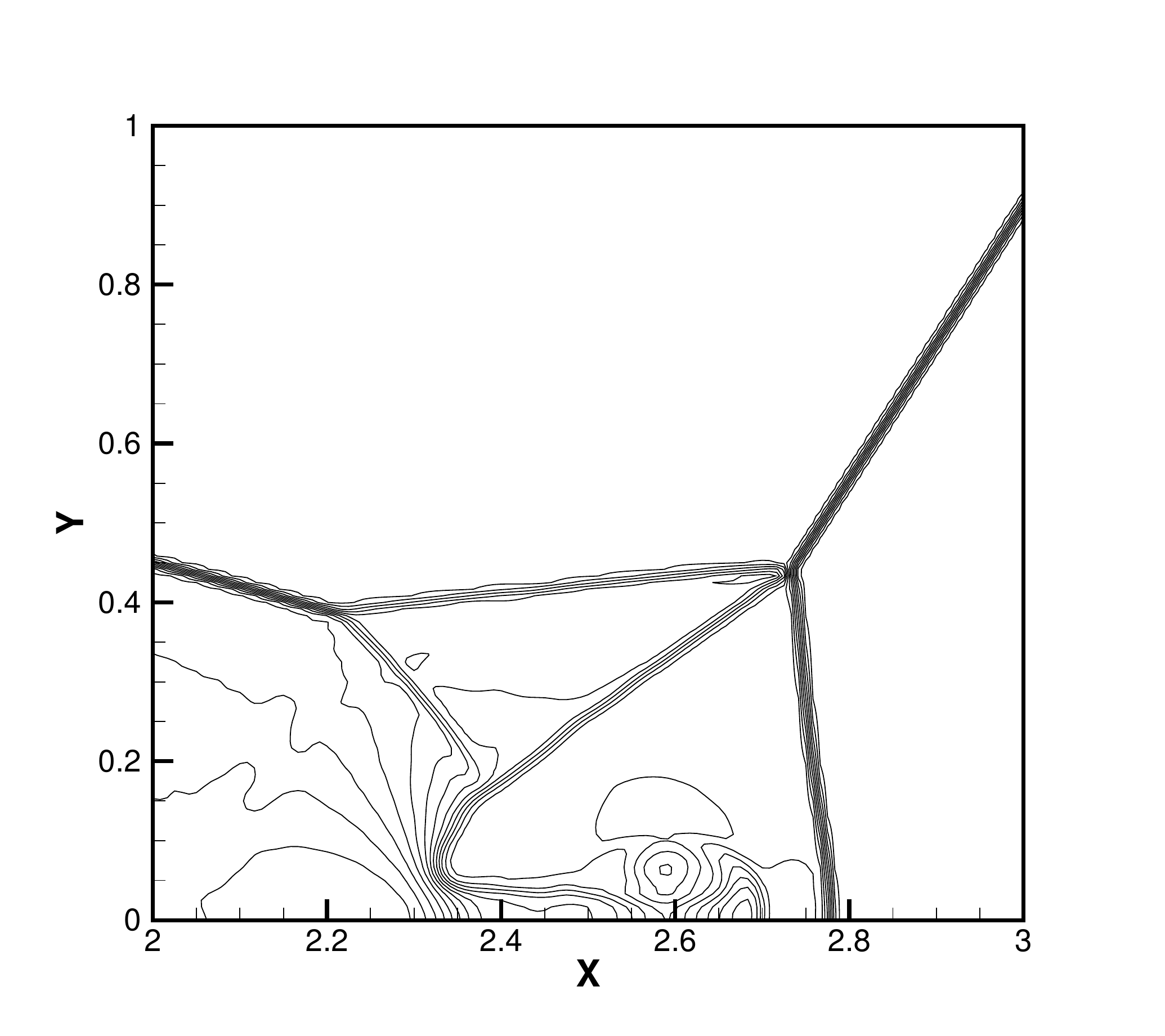}}\quad
 \subfigure[RK-HWENO, $N=480\times 120$]
 {\includegraphics[width=8cm]{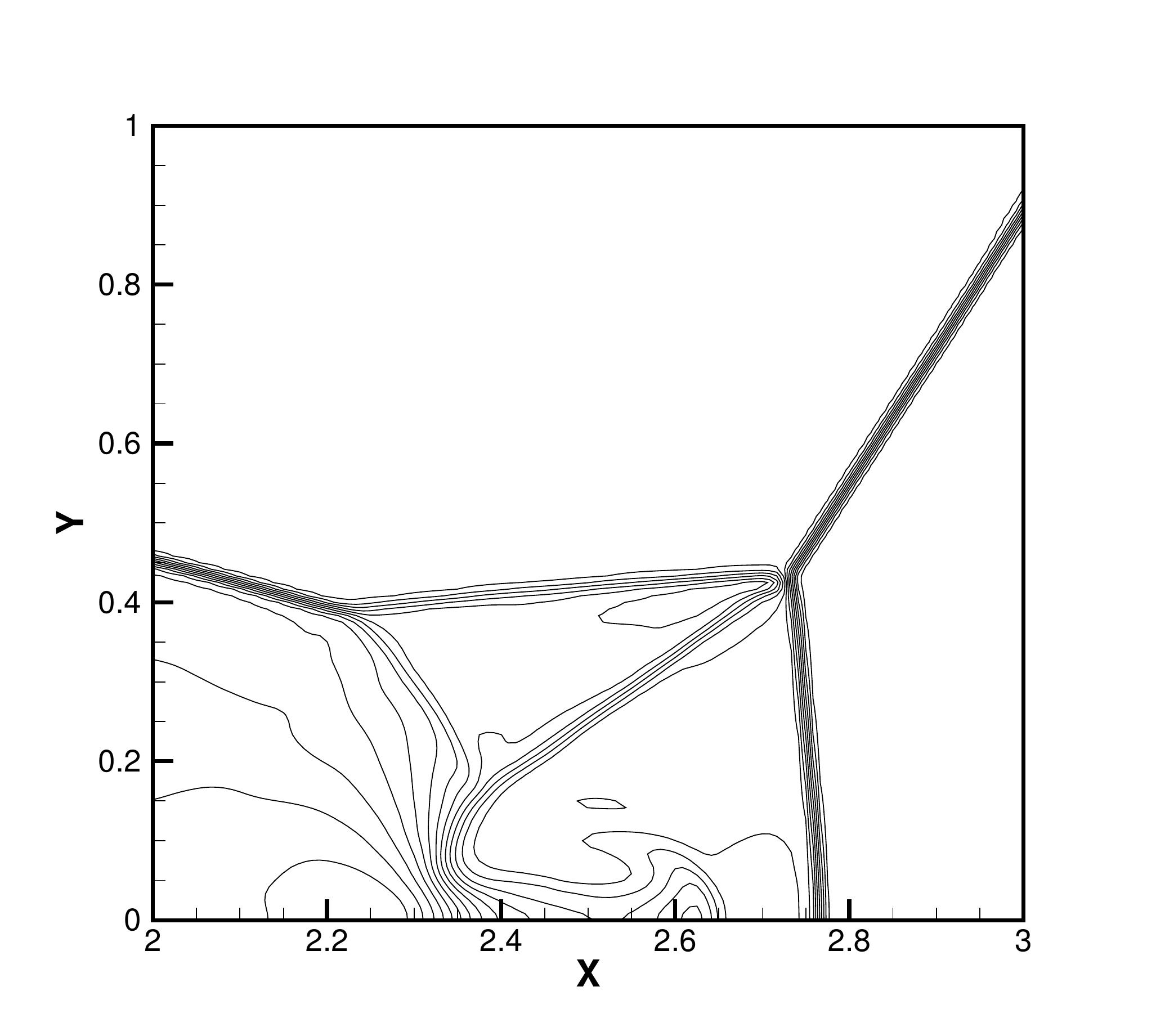}}
}
 \mbox{ \subfigure[ADER-SHWENO,  $N=960\times 240$]
{\includegraphics[width=8cm]{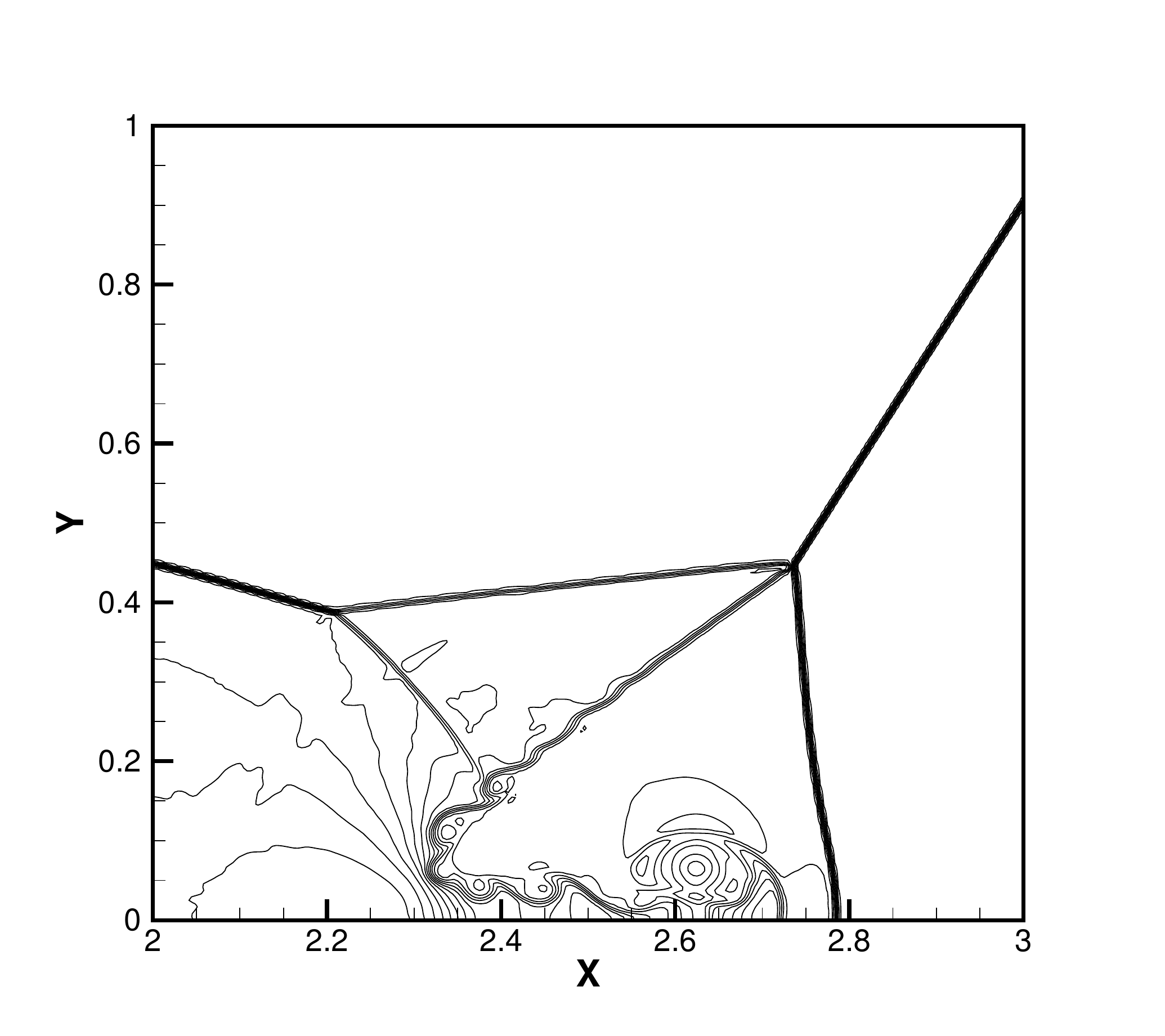}}\quad
 \subfigure[RK-HWENO,  $N=960\times 240$]
{\includegraphics[width=8cm]{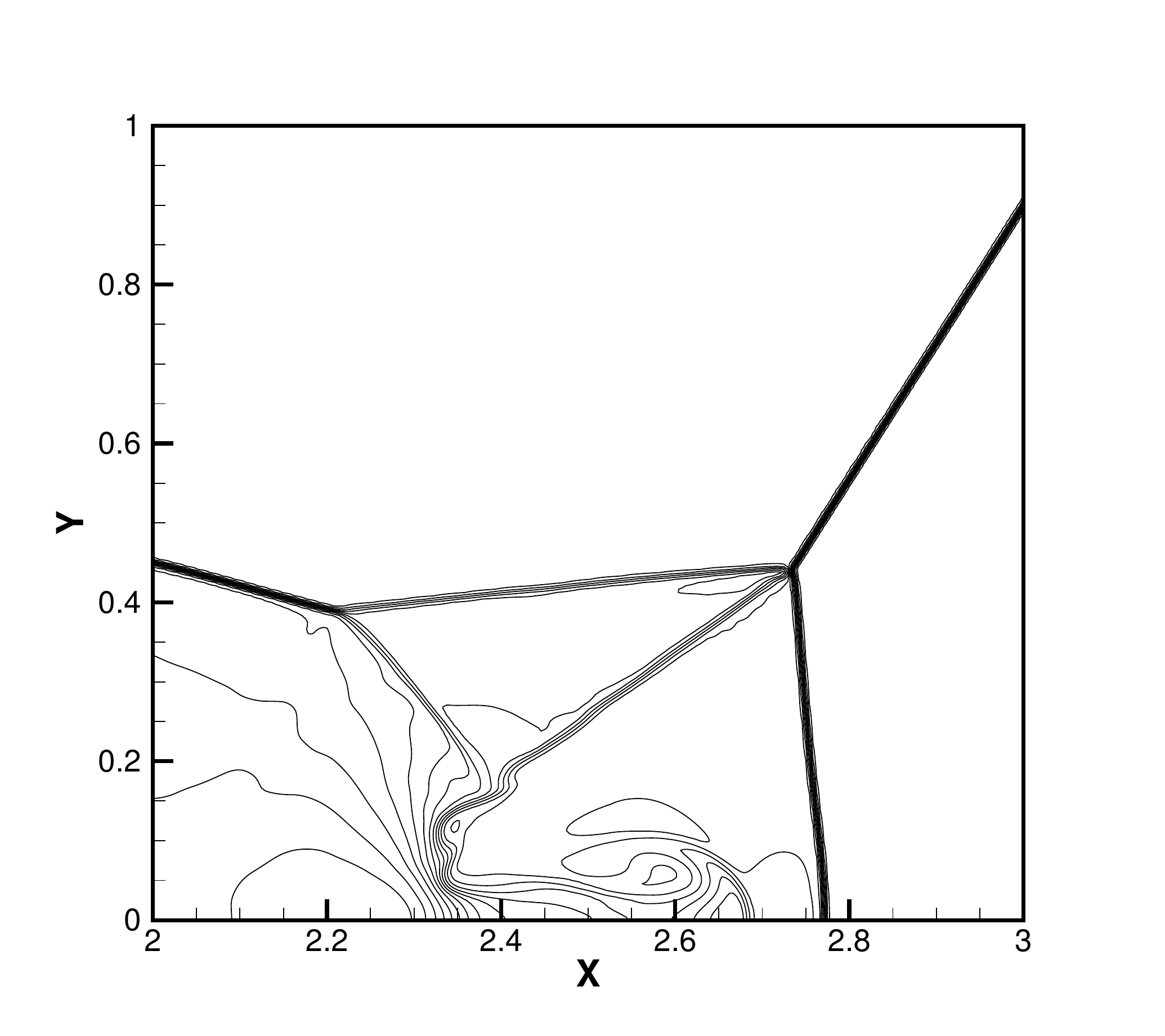}}
}

   \caption{Example~\ref{doublemach} 30 density contours from $2$ to 22. Zoom of Fig. \ref{figdm} in the complex region.}
   \label{figdmzoom}
   \end{center}
   \end{figure}

}
\end{exam}

\section{Conclusions}
\label{secconclusion}
\setcounter{equation}{0}
\setcounter{figure}{0}
\setcounter{table}{0}
 
We have presented a compact and high order ADER scheme using the simple HWENO method for hyperbolic conservation laws. The Lax-Wendroff procedure is adopted in the newly-developed method to convert time derivatives to spatial derivatives, which provides the time evolution of the variables at the cell interface. This information is necessary for the SHWENO reconstructions, which take advantages of the simple WENO \cite{zhu2016,zhu2017} and the classic HWENO \cite{qiu2003,qiu2005}. Compared with the original RK-HWENO \cite{qiu2003,qiu2005}, the new method has two advantages. Firstly, RK-HWENO method must solve the additional equations for reconstructions, which is avoided for the new method. Secondly, the SHWENO reconstruction is performed once with one stencil and is different from the classic HWENO methods, in which both the function and its derivative values are reconstructed with two different stencils, respectively. Thus the new method is more efficient than the RK-HWENO method.  Moreover, the new method makes the best use of the information in the ADER method. Therefore the time evolutions of the cell averages of the derivatives are simpler than that developed in the work \cite{li2021}, where it includes the non-equilibrium and equilibrium parts. Besides, the new method is more compact than the ADER-WENO \cite{titarev2002}. Numerical results in one and two dimensions shown demonstrate the high order for smooth solutions both in space and time and keep non-oscillatory at discontinuities. We recall that the structured mesh has been used in this work. Extending the developed method to the unstructured mesh has  been underway. In addition, we will combine the ADER-SHWENO with the moving mesh method \cite{luo2019} to further improve the resolution.


\begin{thebibliography}{99} 
\bibliographystyle{plain}
\bibliography{reference}

\bibitem{boscheri2018}
W. Boscheri, M. Dumbser, R. Loub\'ere, P.-H. Maire, A second-order cell-centered Lagrangian ADER-MOOD finite volume scheme on multidimensional unstructured meshes for hydrodynamics, {\em Journal of Computational Physics} 358 (2018), 103-129.



\bibitem{castro2011}
M. Castro, B. Costa, W. S. Don, High order weighted essentially non-oscillatory WENO-Z schemes for hyperbolic conservation laws, {\em Journal of Computational Physics} 230 (2011), 1766-1792.




\bibitem{cockburn1989jcp}
B. Cockburn, S.-Y. Lin, C.-W. Shu,
TVB Runge-Kutta local projection discontinuous Galerkin finite element method for conservation laws III: One dimensional systems,
{\em Journal of Computational Physics} 84 (1989), 90-113.


\bibitem{cockburn1998}
B. Cockburn, C.-W. Shu,
The Runge-Kutta discontinuous Galerkin method for conservation laws V: Multidimensional systems,
{\em Journal of Computational Physics} 141 (1998), 199-224.

\bibitem{dumbser2006}
M. Dumbser, Building blocks for arbitrary high order discontinuous Galerkin schemes, {\em Journal of Scientific Computing} 27 (2006), 215-230.

\bibitem{dumbser20071}
M. Dumbser, M. K\"aser, Arbitrary high order non-oscillatory finite volume schemes on unstructured meshes for linear hyperbolic systems, {\em Journal of Computational Physics} 221 (2007), 693-723.

\bibitem{dumbser20072}
M. Dumbser, M. K\"aser, V. A. Titarev, E. F. Toro, Quadrature-free non-oscillatory finite volume schemes on unstructured meshes for nonlinear hyperbolic systems, {\em Journal of Computational Physics} 226 (2007), 204-243.

\bibitem{fambri2017}
F. Fambri, M. Dumbser, O. Zanotti, Space-time adaptive ADER-DG schemes for dissipative flows: Compressible Navier-Stokes and resistive MHD equations, {\em Computer Physics Communications} 220 (2017), 297-318.


\bibitem{gu2017}
Y. Gu, G. Hu, A third order adaptive ADER scheme for one dimensional conservation laws, {\em Communications in Computational Physcis} 22 (2017), 829-851.

\bibitem{harten1987}
A. Harten, B. Engquist, S. Osher, S. R. Chakravarthy, Uniformly high order accurate essentially non-oscillatory schemes III, {\em Journal of Computational Physics} 71 (1987), 231-303.



\bibitem{jiang1996}
G. S. Jiang, C.-W Shu, Efficient implementation of weighted ENO schemes, {\em Journal of Computational Physics} 126 (1996), 202-228.
%

\bibitem{lax1998}
P. D. Lax, X.-D. Liu, Solutions of two-dimensional Riemann problems of gas dynamics by positive schemes, {\em SIAM Journal on Scientific Computing} 19 (1998), 319-340.

\bibitem{levy1999}
D. Levy, G. Puppo, G. Russo, Central WENO schemes for hyperbolic systems of conservation laws, {\em ESAIM: Mathematical Modelling and Numerical Analysis} 33 (1999), 547-571.

\bibitem{levy2000}
D. Levy, G. Puppo, G. Russo, Compact central WENO schemes for multidimensional conservation laws, {\em SIAM Journal on Scientific Computing} 22 (2000), 656-672.

\bibitem{li2016}
J. Li, Z. Du, A two-stage fourth order time-accurate discretization for Lax-Wendroff type flow solvers I. hyperbolic conservation laws, {\em SIAM Journal on Scientific Computing} 38 (2016), A3046-A3069.

\bibitem{li2020}
S. Li, Y. Chen, S. Jiang, An efficient high-order gas-kinetic scheme (I): Euler equations, {\em Journal of Computational Physics} 415 (2020), 109488.

\bibitem{li2021}
S. Li, D. Luo, J. Qiu, Y. Chen, A compact and efficient high-order gas-kinetic scheme, {\em Journal of Computational Physics} 447 (2021), 110661.

\bibitem{liu2014}
N. Liu, H. Tang, A high-order accurate gas-kinetic scheme for one- and two-dimensional flow simulation, {\em Communications in Computational Physics} 15 (2014), 911-943.

\bibitem{luo2019}
D. Luo, W. Huang, J. Qiu, A quasi-Lagrangian moving mesh discontinuous Galerkin method for hyperbolic conservation laws, {\em Journal of Computational Physics} 396 (2019), 544-578.

\bibitem{luo2016}
D. Luo, W. Huang, J.  Qiu, A hybrid LDG-HWENO scheme for KdV-type equations,  {\em Journal of Computational Physics} 313 (2016), 754-774.

\bibitem{pan2016}
L. Pan, K. Xu, Q. Li, J. Li, An efficient and accurate two-stage fourth-order gas-kinetic scheme for the Euler and Navier-Stokes equations, {\em Journal of Computational Physics} 326 (2016), 197-221.




%

%
%
%
%

\bibitem{qiu2003lw}
J. Qiu, C.-W. Shu, Finite difference WENO schemes with Lax-Wendroff-type time discretizations, {\em SIAM Journal on Scientific Computing} 24 (2003), 2185-2198.

\bibitem{qiu2003}
J. Qiu, C.-W. Shu,
Hermite WENO schemes and their application as limiters for Runge-Kutta discontinous Galerkin method: one dimensional case,
{\em Journal of Computational Physics} 193 (2004), 115-135.

\bibitem{qiu2005}
J. Qiu, C.-W. Shu,
Hermite WENO schemes and their application as limiters for Runge-Kutta discontinous Galerkin method II: two dimensional case,
{\em Computers } \& {\em Fluids} 34 (2005), 642-663.

\bibitem{qiu2007lw}
J. Qiu, Hermite WENO schemes with Lax-Wendroff type time discretizations for Hamilton-Jacobi equations, {\em Journal of Computational Mathematics} 25 (2007), 131-144.

\bibitem{EH13}
J. Qiu, C.-W. Shu,
Runge-Kutta discontinuous Galerkin method using WENO limiters,
{\em SIAM Journal on Scientific Computing} 26 (2005), 907-929.

\bibitem{qiu2005lw}
J. Qiu, M. Dumbser, C.-W. Shu, The discontinuous Galerkin method with Lax-Wendroff type time discretizations, {\em Computer Methods in Applied Mechanics and Engineering} 194 (2005), 4528-4543.

\bibitem{ren2015}
X. Ren, K. Xu, W. Shyy, C. Gu, A multi-dimensional high order discontinuous Galerkin method based on gas kinetic theory for viscous flow computations,  {\em Journal of Computational Physics} 292 (2015), 176-193.


%
%
%
%

\bibitem{shu1988}
C.-W. Shu, Total-variation-diminishing time discretizations, {\em SIAM J. Sci. Stat. Comput.} 9 (1988), 1073-1084.

\bibitem{titarev2002}
V. A. Titarev, E. F. Toro, ADER: Arbitrary high order Godunov approach, {\em Journal of Scientific Computing} 17 (2002), 609-618.

\bibitem{titarev2005}
V. A. Titarev, E. F. Toro, ADER schemes for three-dimensional non--linear hyperbolic systems, {\em Journal of Computational Physics} 204 (2005), 715-736

\bibitem{toro2009}
E. F. Toro, Riemann solvers and numerical methods for fluid dynamics: A practical introduction, Springer-Verlag Berlin Heidelberg, 2009.

\bibitem{toro2006}
E. F. Toro, V. A. Titarev, Derivative Riemann solvers for systems of conservation laws and ADER methods, {\em Journal of Computational Physics} 212 (2006), 150-165.

\bibitem{toro2001}
E. F. Toro, R. C. Millington, L. A. M. Nejad, Towards very high-order Godunov schemes. In {\em Godunov Methods: Theory and Applications. Edited Review,}, Kluwer Academic/Plenumm Publishers, 2001, 905-937.

\bibitem{woodward1984}
P. R. Woodward, P. Colella, The numerical simulation of two-dimensional fluid flow with strong shocks, {\em Journal of Computational Physics} 54 (1984), 115-173.

\bibitem{zhao2020}
Z. Zhao, Y. Chen, J. Qiu, A hybrid Hermite WENO scheme for hyperbolic conservation laws, {\em Journal of Computational Physics} 405 (2020), 109175.

\bibitem{zhu2016}
J. Zhu, J. Qiu, A new fifth order finite difference WENO schemes for solving hyperbolic conservation laws, {\em Journal of Computational Physics} 318 (2016), 110-121.

\bibitem{zhu2017}
J. Zhu, J. Qiu, A new type of finite volume WENO schemes for hyperbolic conservation laws, {\em Journal of Scientific Computing} 73 (2017), 1338-1359.

\bibitem{zhu2009}
J. Zhu, J. Qiu, Hermite WENO schemes and their application as limiters for Runge-Kutta discontinuous Galerkin method, III: unstructured meshes, {\em Journal of Scientific Computing} 39 (2009), 293-321.



%
%
%
%





%
%




%
%
%

%
%
%

\end{thebibliography}
\end{document}